\documentclass[11pt]{article}
\usepackage{amssymb,amsmath,amsthm}
\newcommand{\EE}{\mathcal{E}}
\newcommand{\FF}{\mathcal{F}}
\newcommand{\OO}{\mathcal{O}}
\newcommand{\VV}{\mathcal{V}}
\newcommand{\TT}{\mathcal{T}}
\newcommand{\C}{\mathbb{C}}
\newcommand{\N}{\mathbb{N}}
\newcommand{\R}{\mathbb{R}}
\newcommand{\T}{\mathbb{T}}
\newcommand{\Z}{\mathbb{Z}}
\renewcommand{\S}{\mathbb{S}}
\newcommand{\dd}{\,{\rm d}}

\newcommand{\dist}{{\mathrm{dist}}}
\renewcommand{\div}{{\mathrm{div}}}
\newcommand{\curl}{\mathop{\mathrm{curl}}}

\newcommand{\essup}{\mathop{\mathrm{ess\,sup}}}
\renewcommand{\:}{\thinspace :}
\newcommand{\meas}{\mathrm{meas}}
\newcommand{\loc}{\mathrm{loc}}
\newcommand{\ul}{\mathrm{ul}}
\newcommand{\bu}{\mathrm{bu}}
\newcommand{\DS}{\displaystyle}
\newcommand{\1}{\mathbf{1}}

\renewcommand{\Re}{\mathop{\mathrm{Re}}}

\newtheorem{thm}{Theorem}[section]
\newtheorem{df}[thm]{Definition}
\newtheorem{prop}[thm]{Proposition}
\newtheorem{lem}[thm]{Lemma}
\newtheorem{cor}[thm]{Corollary}
\newtheorem{conj}[thm]{Conjecture}
\theoremstyle{definition}
\newtheorem{rem}[thm]{Remark}
\newtheorem{rems}[thm]{Remarks}
\newtheorem{ex}[thm]{Example}
\newtheorem{exs}[thm]{Examples}
\newcommand{\QED}{\mbox{}\hfill$\Box$}

\begin{document}

\title{Distribution of Energy and Convergence to Equilibria 
in Extended Dissipative Systems}

\author{
\null\\
{\bf Thierry Gallay}\\ 
Institut Fourier\\
Universit\'e de Grenoble \& CNRS\\
BP 74\\
38402 Saint-Martin-d'H\`eres, France\\
{\tt Thierry.Gallay@ujf-grenoble.fr}
\and
\\
{\bf Sini\v{s}a Slijep\v{c}evi\'c}\\    
Department of Mathematics\\ 
University of Zagreb\\
Bijeni\v{c}ka 30\\
10000 Zagreb, Croatia\\
{\tt slijepce@math.hr}}

\date{December 4, 2012}

\maketitle

\begin{abstract}
We are interested in understanding the dynamics of dissipative partial
differential equations on unbounded spatial domains.  We consider
systems for which the energy density $e \ge 0$ satisfies an evolution
law of the form $\partial_t e = \div_x f - d$, where $-f$ is the
energy flux and $d \ge 0$ the energy dissipation rate. We also suppose
that $|f|^2 \le b(e)d$ for some nonnegative function $b$. Under these
assumptions we establish simple and universal bounds on the
time-integrated energy flux, which in turn allow us to estimate the
amount of energy that is dissipated in a given domain over a long
interval of time. In low space dimensions $N \le 2$, we deduce that
any relatively compact trajectory converges on average to the set of
equilibria, in a sense that we quantify precisely.  As an application,
we consider the incompressible Navier-Stokes equation in the infinite
cylinder $\R \times \T$, and for solutions that are merely bounded we
prove that the vorticity converges uniformly to zero on large
subdomains, if we disregard a small subset of the time interval.
\end{abstract}

\section{Introduction}\label{intro}
Many time-dependent partial differential equations arising in
Mathematical Physics are dissipative in the sense that there exists a
nonnegative {\em energy density} $e(x,t)$, depending on the space
variable $x \in \R^N$ and the time $t$, which is locally dissipated
under the evolution defined by the system. By this we mean that
$e(x,t)$ satisfies an equation of the form
\begin{equation}\label{eb} 
  \partial_t e(x,t) \,=\, \div_x f(x,t) - d(x,t)~,
\end{equation}
where $-f(x,t) \in \R^N$ denotes the {\em energy flux} in 
the system and $d(x,t) \ge 0$ is the {\em energy dissipation rate}. 
Equivalently, integrating \eqref{eb} with respect to both variables
$x,t$ and applying the divergence theorem, we obtain the 
energy balance equation
\begin{equation}\label{EB} 
  \int_\Omega e(x,T_2)\dd x - \int_\Omega e(x,T_1)\dd x \,=\, 
  \int_{T_1}^{T_2} \int_{\partial\Omega} f(x,t)\cdot\nu \dd\sigma
  \dd t -  \int_{T_1}^{T_2} \int_\Omega d(x,t) \dd x \dd t~,
\end{equation}
which holds for all $T_1 < T_2$ and all admissible domain $\Omega 
\subset \R^N$. Here $\nu$ denotes the outward pointing unit normal on 
$\partial\Omega$, and $\dd\sigma$ is the elementary surface area. 

As a typical example, consider the reaction-diffusion equation 
\begin{equation}\label{rd} 
  \partial_t u(x,t) \,=\, \Delta_x u(x,t) - V'(u(x,t))~, 
  \qquad x \in \R^N~, \quad t > 0~,
\end{equation}
where $u : \R^N \times \R_+ \to \R$ and $V : \R \to \R_+$ is a smooth 
potential. This equation appears for instance in the theory of phase 
transitions \cite{AC} and in population genetics \cite{AW}. In the 
particular case where $V(u) = \frac14 (1-u^2)^2$, Eq.~\eqref{rd} 
is often referred to as the Allen-Cahn equation or the real Ginzburg-Landau
equation. If $u(x,t)$ is any smooth solution of \eqref{rd}, we define
the energy density, the (backward) energy flux, and the energy dissipation 
rate by the formulas
\begin{equation}\label{efdrd} 
  e\,=\, \frac12 |\nabla u|^2 + V(u)~, \qquad f \,=\, u_t \nabla u~,
  \qquad d \,=\, u_t^2~,
\end{equation}
where $u_t = \partial_t u$. It is then straightforward to verify 
that \eqref{eb} holds, which means that energy is locally dissipated
under the evolution defined by \eqref{rd}. Since $V$ is nonnegative, 
we also deduce from \eqref{efdrd} that $|f|^2 \le 2 e d$. We shall 
list in Section~\ref{secEDS} several other examples of classical PDE's 
which define dissipative dynamical systems in the same sense. In most 
of these examples, the energy flux happens to satisfy an inequality of 
the form
\begin{equation}\label{fbound}
  |f|^2 \,\le\, C e d~,
\end{equation}
for some positive constant $C$. 

If a dissipative PDE such as \eqref{rd} is considered in a bounded
domain $\Omega \subset \R^N$, with boundary conditions ensuring 
that $f \cdot \nu \le 0$ on $\partial\Omega$, then \eqref{EB} 
shows that the total energy
\[
  E(t) \,=\, \int_\Omega e(x,t)\dd x 
\]
is a {\em Lyapunov function} of the system, namely $E(t)$ is a 
decreasing function of time for all solutions of \eqref{rd} which 
are not equilibria. Under natural coercivity assumptions on the 
potential $V$, this gradient structure implies that all finite-energy 
solutions of \eqref{rd} in a bounded domain $\Omega$ converge to the 
set of equilibria as $t \to \infty$ \cite{Ha,He}. 

The situation is rather different if we work in an unbounded domain
such as $\Omega = \R^N$. In that case, Eq.~\eqref{rd} may have
travelling wave solutions of the form $u(x,t) = v(x-ct)$ for some
nonzero $c \in \R^N$ \cite{AW}, and such solutions do not converge
uniformly to equilibria as $t \to \infty$. One may object that
travelling waves do converge to equilibria uniformly on compact sets,
but it is possible to construct more complicated solutions for which
convergence to equilibria does not hold even in that weaker sense, see
Example~\ref{ERex} below. Moreover, if $N > 2$, one can exhibit scalar
reaction-diffusion equations of the form $\partial_t u = \Delta u +
F(x,u)$ which have nontrivial time-periodic solutions \cite{GS}. This
is in sharp contrast with what happens for gradient systems, but one
should keep in mind that all counter-examples above involve
infinite-energy solutions.

When the total energy cannot be used as a Lyapunov function, a natural
idea is to exploit the energy balance equation \eqref{EB} or its
differential version \eqref{eb} to obtain relevant information on the 
dynamics. In the context of extended dissipative systems, this 
approach was initiated in a previous paper by the authors \cite{GS}, 
the main conclusions of which can be summarized as follows. If 
$N \le 2$, the reaction-diffusion equation \eqref{rd} on $\R^N$ 
cannot have any nontrivial solution such that $e(x,T_2) \ge e(x,T_1)$ 
for some $T_2 > T_1$ and all $x \in \R^N$; in particular, nontrivial
time-periodic solutions are excluded. Furthermore, all bounded 
solutions converge on average (in time), uniformly on compact sets
(in space), toward the set of equilibria as $t \to +\infty$. In 
other words, due to the local energy dissipation law \eqref{eb}, 
Eq.~\eqref{rd} retains some dynamical properties of usual gradient 
systems, provided $N \le 2$. In contrast, if $N > 2$, highly non-gradient 
behaviors such as nontrivial periodic orbits can occur. The conclusions 
of \cite{GS} also apply to a damped hyperbolic equation which reduces to 
\eqref{rd} in the limit of strong damping. 

The aim of the present paper is to strengthen and generalize the
results of \cite{GS}. Instead of considering a particular equation, we
work in an abstract setting, assuming only the local energy
dissipation law \eqref{eb} and an estimate of the form \eqref{fbound}
for the energy flux. As a consequence, our results apply to a much
larger class of systems, some of which are listed in
Section~\ref{secEDS}. Another substantial progress with respect to
\cite{GS} is a new estimate on the time-integrated energy flux through
a closed hypersurface, which we derive in Section~\ref{secFlux}. This
bound allows us to obtain {\em quantitative} versions of the main
results in \cite{GS}. For instance, in Section~\ref{secEn} we give an
explicit estimate of the energy dissipated in a given domain over a
long time interval, and in Section~\ref{secConv} we measure the
fraction of time spent by any relatively compact trajectory outside
a neighborhood of the set of equilibria. As can be expected from 
\cite{GS}, our results depend strongly on the space dimension $N$, and 
some of them even fail if $N > 2$. As a final application, we consider in
Section~\ref{secNS} the incompressible Navier-Stokes equation in the
two-dimensional cylinder $\R \times \T$, and for solutions that are
merely bounded we prove some convergence results for the vorticity
which are apparently new in this context. 

\section{Extended Dissipative Systems}\label{secEDS}

To treat in a unified way various dissipative PDE's on unbounded 
domains, we introduce in this section the notion of an {\em 
extended dissipative system}, which will be studied in the rest
of this paper. We also list a few classical examples which fit
into our abstract framework. 

Let $X$ be a metrizable topological space. We say that a family
$(\Phi(t))_{t\ge0}$ of continuous maps in $X$ is a {\em continuous
semiflow} on $X$ if\\[1mm]
\null\hskip 8pt i) $\Phi(0) = \1$ (the identity map);\\[1mm]
\null\hskip 5pt ii) $\Phi(t_1+t_2) = \Phi(t_1)\circ\Phi(t_2)$ 
for all $t_1,t_2 \ge 0$;\\[1mm]
\null\hskip 2pt iii) For any $T > 0$, the map $(t,u) \mapsto 
\Phi(t)u$ is continuous from $[0,T]\times X$ to $X$.

In particular, if $u_0 \in X$, the trajectory $u : \R_+ \to X$ defined
by $u(t) = \Phi(t)u_0$ for all $t \in \R_+ = [0,\infty)$ is
continuous, and $u(t)$ depends continuously on the initial data $u_0$,
uniformly in time on compact intervals. As an example, if $V(u) =
\frac14(1-u^2)^2$, the reaction-diffusion equation \eqref{rd} defines
a continuous evolution semiflow on the space $X = C^0(\R^N) \cap
L^\infty(\R^N)$, if $X$ is equipped with the topology of uniform
convergence on compact sets of $\R^N$. The systems we are interested
in are those for which one can define an energy density $e$, an energy
flux $f$, and an energy dissipation rate $d$ with the same properties
as in the example \eqref{rd}. This leads to the following definition\:

\begin{df}\label{edsdef}
Let $N \in \N^*$. We say that a continuous semiflow $(\Phi(t))_{t\ge0}$
on a metrizable space $X$ is an {\em extended dissipative system} on 
$\R^N$ if one can associate to each $u \in X$ a triple $(e,f,d)$ 
with $e,d \in C^0(\R^N,\R_+)$ and $f \in C^0(\R^N,\R^N)$ such 
that\:\\[1mm]
\null\hskip 5pt A1\: The functions $e,f,d$ depend continuously on 
$u \in X$, uniformly on compact sets of $\R^N$;\\[1mm]
\null\hskip 5pt A2\: $\,|f|^2 \le b(e)d$ for some nondecreasing 
function $b: \R_+ \to \R_+$;\\[1mm]
\null\hskip 5pt A3\: $\,d \equiv 0$ only if $\Phi(t)u = u$ for all 
$t \ge 0$;\\[1mm]
\null\hskip 5pt A4\: Under the evolution defined by the semiflow 
$(\Phi(t))_{t\ge0}$, the time-dependent quantities $e,f,d$ \\
\null\hskip 30pt satisfy the energy balance $\partial_t e = \div f - d$ 
in the sense of distributions on $\R^N \times \R_+$.
\end{df}

\begin{rems}\label{edsrem}\quad\\[1mm]
{\bf 1.} More generally, one can define extended dissipative systems
on any (unbounded) domain $\Omega \subset \R^N$ by substituting 
$\Omega$ for $\R^N$ everywhere in Definition~\ref{edsdef}. For simplicity, 
we restrict ourselves to the case of the whole space $\R^N$ in this 
paper, except in Section~\ref{secNS} where we shall consider a system 
defined in a cylindrical domain.

\noindent{\bf 2.} We emphasize that, in Definition~\ref{edsdef}, both
the energy density $e$ and the energy dissipation rate $d$ are
supposed to be nonnegative. The first condition ensures that the
energy density is bounded from below, and the positivity of $d$
together with the energy balance \eqref{eb} imply that energy is
locally {\em dissipated} (and never created) in the system.

\noindent{\bf 3.} In assumption~A2, it is understood that the function
$b : \R_+ \to \R_+$ is independent of $u \in X$. In many examples one
can take $b(e) = Ce$ for some positive constant $C$, as in
\eqref{fbound}, but the generalization proposed here is necessary if
one considers systems such as the nonlinear diffusion equation (see
below) or the two-dimensional vorticity equation (see
Section~\ref{secNS}).

\noindent{\bf 4.} Assumption~A3 means that all trajectories, except
equilibria, dissipate some energy. Note that, if $d \equiv 0$, then $f
\equiv 0$ by A2, hence \eqref{eb} already implies that the energy
density is time-independent. We emphasize that an extended dissipative
system may have equilibria for which $d$ is not identically zero;
these may be called {\em nonequilibrium steady states}, in the
terminology of Statistical Mechanics. On the other hand, if one
considers systems with a continuous group of symmetries, it may be
useful to relax assumption A3 so as to allow for a vanishing energy
dissipation on {\em relative equilibria} of the system; these are
equilibria up to the action of the symmetry group, see the example of
the complex Ginzburg-Landau equation below. 

\noindent{\bf 5.} Our final comment on Definition~\ref{edsdef}
concerns the regularity of $e$, $f$, and $d$. To avoid technicalities,
we have supposed that, for each $u \in X$, the densities $e,f,d$ are
continuous functions on $\R^N$. Moreover, if $u(t)$ varies
continuously in time, the associated quantities $e(t),f(t),d(t)$ are
jointly continuous in space and time, in view of assumption~A1. In the
particular case where $u(t) = \Phi(t)u_0$ for some $u_0 \in X$,
assumption~A4 implies that the integrated energy balance equation
\eqref{EB} holds for all $T_2 > T_1 \ge 0$ and any smooth bounded
domain $\Omega \subset \R^N$, and since $\Phi$ is a continuous
semiflow we even know that all four terms in \eqref{EB} depend
continuously on the initial data $u_0 \in X$. These comfortable
assumptions are not unrealistic, and can be verified in all systems
listed below if we choose functions spaces of sufficiently high
regularity. However, especially in nonparabolic PDE's, it is often
more convenient to use larger function spaces, in which (for instance)
the energy density is locally integrable but not continuous. In that
case, intead of A1 and A4, it is sufficient to require that the energy
balance equation \eqref{EB} be satisfied for all $T_2 > T_1 \ge 0$ and
any smooth bounded domain $\Omega \subset\R^N$, and that the various
quantities in \eqref{EB} depend continuously on the initial data.
\end{rems}

In complement to Definition~\ref{edsdef}, we specify two properties
of extended dissipative systems which will play an important role
in Sections~\ref{secFlux} to \ref{secConv}.

\begin{df}\label{bcdef}
We say that an extended dissipative system in the sense of
Definition~\ref{edsdef} is\\[1mm]
\null\hskip 5pt $\bullet$ {\em bounded}, if there exists $E > 0$ 
such that $e(x) \le E$ for all $x \in \R^N$ and all $u \in X$;\\[1mm]
\null\hskip 5pt $\bullet$ {\em compact}, if the space $X$ is compact. 
\end{df}

\begin{rem}\label{bdrem}
In many extended dissipative systems, boundedness can be achieved
simply by restricting the space $X$ to a subset which is positively
invariant under the evolution defined by the semiflow, and on which
the energy density is uniformly bounded by some positive constant
$E$. In that case, we can replace $b(e)$ by a positive constant $\beta
\ge b(E)$ in assumption~A2.  We thus obtain the relation $d \ge
|f|^2/\beta$, which quantifies how much energy is dissipated in the
system when the energy flux is not identically zero. On the other
hand, compactness can often be achieved by endowing $X$ with a
sufficiently weak topology. Note however that, with the definitions
above, a compact extended dissipative system is not necessarily
bounded.
\end{rem}

To substantiate Definitions~\ref{edsdef} and \ref{bcdef}, we
now give a few concrete examples. 

\begin{exs}\label{edsex}\quad\\
{\bf 1.} {\em A reaction-diffusion equation}\\
We consider again the reaction-diffusion equation \eqref{rd}, and
specify in which function spaces it defines an extended dissipative
system in the sense of Definition~\ref{edsdef}. There are of course
many possibilities, and we just mention here two reasonable ones.
Since we want global solutions of \eqref{rd}, it is natural to assume
that the potential $V : \R \to \R_+$ is coercive in some sense. For
instance, we can suppose that there exists a constant $m \ge 0$ such
that $u V'(u) \ge m$ whenever $|u|$ is sufficiently large. Then it is
known that the Cauchy problem for \eqref{rd} is globally well-posed in
$C^k_\bu(\R^N)$, the Banach space of all functions $u : \R^N \to \R$
that are bounded and uniformly continuous together with their
derivatives up to order $k \in \N$. This means that \eqref{rd} defines
a continuous semiflow $(\Phi(t))_{t\ge0}$ on $X =
C^k_\bu(\R^N)$. Moreover, if $k \ge 2$, the quantities $e,f,d$
introduced in \eqref{efdrd} belong to $C^0(\R^N)$ and depend
continuously on $u \in X$, uniformly on compact sets of
$\R^N$. Together with \eqref{eb} and \eqref{fbound}, this shows that,
if $k \ge 2$, the semiflow of \eqref{rd} on $X = C^k_\bu(\R^N)$ is an
extended dissipative system. 

This system is not bounded in the sense of Definition~\ref{bcdef}, but
it becomes so if we choose for $X$ any bounded subset $B \subset
C^k_\bu(\R^N)$ that is positively invariant under the semiflow.  In
addition, if $k > 2$ and if we consider $B$ as a subset of
$C^0(\R^N)$, equipped with the topology of uniform convergence on
compact sets, then the closure of $B$ in $C^0(\R^N)$ is compact and
the restriction of the semiflow $(\Phi(t))_{t\ge0}$ to $\bar B$
defines a compact extended dissipative system. The idea of introducing
a localized topology to restore compactness plays an important role in
the study of dissipative PDE's on unbounded domains, in particular
when constructing global attractors \cite{BV,Fe,MS}.

Instead of $C^k_\bu(\R^N)$, another possible choice is the uniformly
local Sobolev space $H^s_\ul(\R^N)$, on which \eqref{rd} also defines
a continuous semiflow if $s > N/2$, see \cite{ABCD,GS}. If moreover $s
> 2+N/2$, the densities \eqref{efdrd} are continuous and we again
obtain an extended dissipative system. As above, if we restrict our
analysis to a bounded invariant subset $B \subset H^s_\ul(\R^N)$ and
if we take the closure of $B$ with respect to the topology of
$L^2_\loc(\R^N)$, the restriction of the semiflow to $\bar B$ defines
an extended dissipative system that is bounded and compact in the
sense of Definition~\ref{bcdef}.

\medskip\noindent{{\bf 2.} {\em A strongly damped wave equation}}
\cite{Ma,PZ}\\
Given $\alpha \ge 0$ and a smooth potential $V : \R \to \R_+$, we 
consider the equation
\begin{equation}\label{dhe}
  u_{tt} + u_t - \alpha \Delta u_t \,=\, \Delta u - V'(u)~, 
  \qquad x \in \R^N~, \quad t > 0~,  
\end{equation}
where $u : \R^N \times \R_+ \to \R$. As usual, this second-order 
equation can be written as a first-order system for the pair 
$(u,u_t)$. For simplicity, we assume that the potential $V$ is 
quadratic near infinity, namely $V''(u) = m > 0$ for all sufficiently
large $u \in \R$. Then the initial value problem for Eq.~\eqref{dhe} is 
globally well-posed in the uniformly local space $X = H^s_\ul(\R^N) 
\times H^{s-1}_\ul(\R^N)$ if $s > N/2$. For any pair $(u,u_t) \in X$, 
we introduce the densities
\[
  e\,=\, \frac12 u_t^2 + \frac12 |\nabla u|^2 + V(u)~, \qquad 
  f \,=\, u_t (\nabla u + \alpha\nabla u_t)~, \qquad d \,=\, 
  u_t^2 + \alpha |\nabla u_t|^2~, 
\]
which are well-defined and continuous provided $s > 2 + N/2$, or $s >
1 +N/2$ if $\alpha = 0$. Then \eqref{dhe} implies that the local
energy dissipation law \eqref{eb} is satisfied, and a direct calculation 
shows that \eqref{fbound} holds with $C = 2\max(1,\alpha)$. Finally 
$d \equiv 0$ implies $u_t \equiv 0$. Thus the semiflow of the
strongly damped wave equation \eqref{dhe} in $X$ is an extended dissipative
system in the sense of Definition~\ref{edsdef}. In the particular
case where $\alpha = 0$, the local dissipation of energy for  
Eq.~\eqref{dhe} was already studied in \cite{GS}.

\medskip\noindent{{\bf 3.} {\em The complex Ginzburg-Landau equation}}
\cite{AK,Co,Mi}\\
Our next example originates from the complex Ginzburg-Landau equation
\begin{equation}\label{cgl}
  u_t \,=\, (1+i\alpha) \Delta u + u - (1+i\beta)|u|^2u~, 
  \qquad x \in \R^N~, \quad t > 0~, 
\end{equation}
where $u : \R^N \times \R_+ \to \C$ and $\alpha,\beta$ are real 
parameters. To have a gradient structure, we assume that $\beta = \alpha$, 
and we introduce the auxiliary function $v(x,t) \,=\, u(x,t) e^{i\alpha t}$, 
which satisfies the equation
\begin{equation}\label{cgl2}
  v_t \,=\, (1+i\alpha) \Bigl(\Delta v + v - |v|^2v\Bigr)~, 
  \qquad x \in \R^N~, \quad t > 0~. 
\end{equation}
The Cauchy problem for \eqref{cgl2} is globally well-posed in the 
function $X = C^k_\bu(\R^N,\C)$ for $k \ge 0$ or $X = H^s_\ul(\R^N,\C)$ 
for $s > N/2$. If in addition $k \ge 2$ or $s > 2+N/2$, then for any 
$v \in X$ the densities
\[
  e\,=\, \frac12 |\nabla v|^2 + \frac14 (1-|v|^2)^2~, 
  \qquad f \,=\, \Re (v_t \nabla \bar v)~, \qquad d \,=\, 
  \frac{|v_t|^2}{1 + \alpha^2}~,
\]
are well-defined and continuous. A direct calculation also shows that
\eqref{eb} holds, as well as \eqref{fbound} with $C = 2(1+\alpha^2)$.
Thus the semiflow of \eqref{cgl2} in $X$ is an extended dissipative
system in the sense of Definition~\ref{edsdef}. We also remark that $d
\equiv 0$ if and only if $u(x,t) = v(x) e^{-i\alpha t}$ for some $v
\in X$, which means that $u(\cdot,t)$ is a {\em relative equilibrium}
of \eqref{cgl}\: $u(\cdot,t)$ moves without dissipation along an orbit
of the symmetry group $U(1)$. Thus the semiflow of \eqref{cgl} in $X$
is an extended dissipative system only if assumption A3 is relaxed as
suggested in Remark~\ref{edsrem}.4. 

\medskip\noindent{{\bf 4.} {\em The Landau-Lifshitz-Gilbert equation}}
\cite{GD}\\ 
We now consider a vector-valued PDE appearing in micromagnetism. Given 
$\alpha\in \R$, the Landau-Lifshitz equation reads
\begin{equation}\label{llg}
  u_t \,=\, - u \wedge (u \wedge \Delta u) + \alpha u \wedge \Delta u~,
  \qquad x \in \R^N~, \quad t > 0~, 
\end{equation}
where $u : \R^N \times \R_+ \to \S^2 = \{v \in \R^3\,|\, |v| =
1\}$. Here $\wedge$ denotes the usual cross product in $\R^3$.  In
particular $- u \wedge (u \wedge \Delta u) = \Delta u - u(u\cdot\Delta
u) = \Delta u + |\nabla u|^2 u$ is the orthogonal projection of
$\Delta u$ onto the plane orthogonal to the direction $u \in \S^2$,
and $u \wedge \Delta u$ is the same vector rotated by $\pi/2$ in the
orthogonal plane. The initial value problem for \eqref{llg} is {\em
  locally} well-posed in the space $X = C^k_\bu(\R^N)$ for $k > 0$ or
$X = H^s_\ul(\R^N)$ for $s > N/2$, but in general finite-time
singularities are expected to occur, unlike in the previous examples. 
To obtain a continuous semiflow $(\Phi(t))_{t\ge0}$, it is therefore 
necessary to restrict our space $X$ to a family of global trajectories. 
Now, if $k \ge 2$ or $s > 2+N/2$, the densities
\[
  e\,=\, \frac12 |\nabla u|^2~, \qquad f \,=\, u_t \nabla u 
  \,\equiv\, \sum_{k=1}^N (\partial_t u_k) \nabla u_k~, 
  \qquad d \,=\, |u \wedge \Delta u|^2~,
\]
are well-defined and continuous for any $u \in X$, and it is again
straightforward to verify that \eqref{eb} and \eqref{fbound} hold with
$C = 2(1+\alpha^2)$.  Moreover $d \equiv 0$ implies $u \wedge \Delta u
\equiv 0$, hence $u_t \equiv 0$. Thus Eq.~\eqref{llg} also defines an
extended dissipative system in the sense of Definition~\ref{edsdef},
provided we restrict the space $X$ to a suitable family of global
solutions.

\medskip\noindent{{\bf 5.} {\em A nonlinear diffusion equation}}\\ 
To motivate assumption~A2 in Definition~\ref{edsdef}, we also give 
an example where the relation between the energy flux and the 
energy disipation is more complex than in \eqref{fbound}. Given a 
smooth function $a : \R \to (0,\infty)$, we consider the nonlinear 
diffusion equation
\begin{equation}\label{nld}
  u_t \,=\, \div(a(u)\nabla u)~, \qquad x \in \R^N~, \quad t > 0~,
\end{equation}
which is globally well-posed in the space $X = C^k_\bu(\R^N)$ for $k 
\ge 0$ or $X = H^s_\ul(\R^N)$ for $s > N/2$. If $k \ge 1$ or 
$s > 1+N/2$, we denote for all $u \in X$\:
\[
  e\,=\, \frac12 u^2~, \qquad f \,=\, u a(u)\nabla u, \qquad
  d \,=\, a(u)|\nabla u|^2~.
\]
Then \eqref{eb} holds, and since $a(u) > 0$ it is clear that 
$d \equiv 0$ implies $u_t \equiv 0$. Moreover $|f|^2 \le 2 a(u)ed$.
Thus, if we define
\[
  b(e) \,=\, 2 e \sup\{a(u)\,|\, u^2 \le 2e\}~, \qquad e \ge 0~,
\]
then $e \mapsto b(e)$ is increasing and $|f|^2 \le b(e)d$ by construction. 
Thus \eqref{nld} defines an extended dissipative system in $X$ in 
the sense of Definition~\ref{edsdef}. 

\medskip\noindent{{\bf 6.} {\em The two-dimensional vorticity equation}}\\ 
As a final example, we consider the vorticity equation associated to 
the two-dimensional incompressible Navier-Stokes system. In this 
model, the velocity of the fluid, which is denoted by 
$u(x,t) \in \R^2$, satisfies the incompressibility condition 
$\partial_1 u_1 + \partial_2 u_2 = 0$, and the corresponding
vorticity field $\omega = \partial_1 u_2 - \partial_2 u_1$ evolves 
according to the advection-diffusion equation
\begin{equation}\label{nsR2}
  \partial_t\omega  + u\cdot \nabla \omega \,=\, \Delta \omega~,
  \qquad x \in \R^2~, \quad t > 0~.  
\end{equation}
If we define the enstrophy density $e$, the enstrophy flux $f$, and 
the enstrophy dissipation rate $d$ by the formulas
\begin{equation}\label{efdnsR2} 
  e\,=\, \frac12 \omega^2~, \qquad f \,=\, \omega\nabla \omega 
  - \frac12 u\,\omega^2~, \qquad d \,=\, |\nabla \omega|^2~,
\end{equation}
it is easy to verify that \eqref{eb} is satisfied. Moreover, $d \equiv
0$ clearly implies that $\partial_t \omega \equiv 0$.  However, it is
not possible to obtain here an inequality of the form \eqref{fbound},
nor of the generalized form adopted in assumption~A2 of
Definition~\ref{edsdef}. The main difficulty comes from the term
$\frac12 u\,\omega^2$ in the enstrophy flux.  Since the velocity $u$
is not a local function of $\omega$, the bound \eqref{fbound} cannot
hold pointwise unless the contribution of $u$ is absorbed into the
constant $C$ in the right-hand side. This requires a uniform bound on
the $L^\infty$ norm of the velocity field, which is not known to hold
for solutions of the two-dimensional Navier-Stokes equations that are
only bounded in space, see \cite{GMS,ST,Ze}. In addition, the term
$\frac12 u\,\omega^2$ does not contain any derivative of $\omega$,
hence does not necessarily vanish when $d = 0$. This means that
enstrophy can (a priori) be transported without any dissipation,
whereas it is essential in our approach that the energy dissipation be
bounded from below in terms of the energy flux. Surprisingly enough,
we shall show in Section~\ref{secNS} that these difficulties
essentially disappear if we consider the vorticity equation
\eqref{nsR2} in the infinite cylinder $\R \times \T$ instead of the
whole plane $\R^2$. Thus, if we assume periodicity in one space 
direction, the semiflow of \eqref{nsR2} defines a one-dimensional 
extended dissipative system which (nearly) satisfies the assumptions
in Definition~\ref{edsdef}. 
\end{exs}

\begin{rem}\label{restrict}
The above list of examples can certainly be made longer, but all 
extended dissipative systems we are aware of are related somehow
to a parabolic equation involving a second order differential operator.
Higher-order systems, such as the Cahn-Hilliard equation, do not fit 
into our framework since they require a radical modification of the bound 
\eqref{fbound}, which would affect our results in an essential 
way. 
\end{rem}

\section{Bounds on the Energy Flux}\label{secFlux}

We now begin our study of the dynamics of extended dissipative
systems. Given a continuous semiflow $(\Phi(t))_{t\ge0}$ on a
metrizable space $X$ satisfying the assumptions of 
Definition~\ref{edsdef} for some $N \in \N^*$, we consider a
trajectory $u(t) = \Phi(t)u_0$ for which the energy density $e(x,t)$
is {\em uniformly bounded}.  This is always the case if our system is
bounded in the sense of Definition~\ref{bcdef}, and as was observed in
Remark~\ref{bdrem} boundedness can often be achieved by restricting
the space $X$ to a suitable positively invariant subset. If we denote
\begin{equation}\label{e0def}
  e_0 \,=\, \sup_{x \in \R^N} e(x,0) \,<\, \infty~, \qquad
  \beta \,=\, \sup_{x \in \R^N} \sup_{t \ge 0} b(e(x,t)) 
  \,<\, \infty~,
\end{equation}
where $e \mapsto b(e)$ is the nondecreasing function appearing
in Definition~\ref{edsdef}, assumption~A2 then implies
\begin{equation}\label{fbound2}
  |f(x,t)|^2 \,\le\, \beta\,d(x,t)~, \qquad x \in \R^N~, \quad 
  t \ge 0~. 
\end{equation}
Using only \eqref{eb}, \eqref{fbound2}, and the positivity of 
$e(x,t)$ and $d(x,t)$, we shall derive a universal bound on the total 
energy flux through a given hypersurface in $\R^N$ during the time 
interval $[0,T]$. 
 
We first consider the one-dimensional case $N = 1$, where our 
hypersurface is reduced to a single point. Our main result in 
this case is\:

\begin{prop}\label{flux1}
Assume that $N = 1$, and let $u(t) = \Phi(t)u_0$ be a trajectory 
for which the energy density $e(x,t)$ satisfies \eqref{e0def}. Then 
for any $x \in \R$ and any $T > 0$, we have
\begin{equation}\label{F1bd}
  \left|\,\int_0^T f(x,t)\dd t \,\right| \,\le\, \sqrt{\beta T e_0}~.
\end{equation}
\end{prop}

\noindent{\bf Proof.} Without loss of generality, we can assume 
that $\beta > 0$. Given any $T > 0$, we introduce the integrated 
energy flux
\begin{equation}\label{intflux1}
  F_1(x,T) \,=\, \int_0^T f(x,t)\dd t~, \qquad x \in \R~,
\end{equation}
which is a continuous function of $x \in \R$. For any $x_0 \in \R$, we
shall prove that $F_1(x_0,T) \le (\beta T e_0)^{1/2}$.  Together with
the corresponding lower bound $F_1(x_0,T) \ge -(\beta T e_0)^{1/2}$,
which can be established in a similar way, this gives the desired
conclusion.

For any $x > x_0$, the energy balance equation \eqref{EB} 
with $\Omega = (x_0,x)$, $T_1 = 0$, $T_2 = T$ implies
\begin{equation}\label{F1a}
  F_1(x,T) \,=\, F_1(x_0,T) + \int_{x_0}^x \Bigl(e(y,T) - e(y,0)\Bigr)
  \dd y + \int_{x_0}^x \int_0^T d(y,t) \dd t \dd y~.
\end{equation}
Since $e(y,T) \ge 0$ and $e(y,0) \le e_0$, the first integral in the
right-hand side is bounded from below by $-e_0(x-x_0)$. On the 
other hand, using \eqref{fbound2} and the Cauchy-Schwarz inequality, 
we obtain
\[
  \int_0^T d(y,t)\dd t \,\ge\, \frac{1}{\beta}\int_0^T f(y,t)^2 
  \dd t \,\ge\, \frac{1}{\beta T} F_1(y,T)^2~,
\]
for all $y \in [x_0,x]$. Thus \eqref{F1a} implies
\begin{equation}\label{F1b}
  F_1(x,T) \,\ge\, F_1(x_0,T) - e_0(x-x_0) + \frac{1}{\beta T} 
  \int_{x_0}^x F_1(y,T)^2 \dd y~, \qquad x > x_0~.
\end{equation}
We now compare $F_1(x,T)$ to the solution $\tilde F(x)$ of the 
differential equation
\begin{equation}\label{F1c}
  \tilde F'(x) \,=\, -e_0 + \frac{1}{\beta T}\tilde F(x)^2~, 
  \qquad x > x_0~,
\end{equation}
with initial data $\tilde F(x_0) = F_1(x_0,T)$. If $F_1(x_0,T) > 
(\beta T e_0)^{1/2}$, then $\tilde F$ is strictly increasing and blows
up at some finite point $x_1 > x_0$. But \eqref{F1b} then implies that
$F_1(x,T) \ge \tilde F(x)$ for all $x \in (x_0,x_1)$, which leads 
to a contradiction because we know that $F_1(x,T)$ is uniformly
bounded on $[x_0,x_1]$. Thus we must have $F_1(x_0,T) \le 
(\beta T e_0)^{1/2}$. 
\QED

\begin{rem}\label{noper1}
Albeit elementary, Proposition~\ref{flux1} already has interesting
dynamical consequences. For instance, it immediately implies that an
extended dissipative system on $\R$ cannot have any nontrivial
time-periodic orbit with uniformly bounded energy, see \cite{GS}.
Indeed, for such a periodic orbit, the last integral in the right-hand
side of \eqref{F1a} grows linearly in $T$ as $T \to \infty$ (if the
interval $[x_0,x]$ is large enough to include a region where
energy dissipation is nonzero), whereas the first integral is
uniformly bounded by periodicity and the flux terms are $\OO(T^{1/2})$
by \eqref{F1bd}. Thus \eqref{F1a} cannot hold for sufficiently large
times.
\end{rem}

We next investigate the analog of Proposition~\ref{flux1} in 
the higher-dimensional case $N \ge 2$. Here we consider the energy 
flux through the boundary of the ball $B_R = \{x \in \R^N\,|\, |x| 
\le R\}$, for various values of the radius $R$. We recall that the 
Euclidean measure of the sphere $\partial B_R$ is $\omega_N R^{N-1}$, 
where 
\[
  \omega_N \,=\, \frac{2\pi^{N/2}}{\Gamma(N/2)}~, \qquad \hbox{and}
  \quad \Gamma(\alpha) = \int_0^\infty t^{\alpha-1} e^{-t}\dd t~, \quad 
  \alpha > 0~.
\]
Given $R > 0$ and $T > 0$, we thus define the integrated flux
\begin{equation}\label{intfluxN}
  F(R,T) \,=\, \int_0^T \!\!\int_{|x|=R} f(x,t)\cdot 
  \frac{x}{|x|}\dd\sigma \dd t~,
\end{equation}
which represents the total energy entering the ball $B_R$ through 
the boundary over the time interval $[0,T]$ (the energy leaving
the ball is of course counted negatively). 

Before stating our result, we introduce the higher-dimensional
analog of the differential equation \eqref{F1c}, which (after 
suitable normalization) becomes
\begin{equation}\label{EDO}
  H'(r) + \frac{N{-}1}{r} H(r) \,=\, -1 + H(r)^2~, \qquad r > 0~.
\end{equation}
The following elementary result will be established in 
Section~\ref{secApp}.

\begin{lem}\label{EDOlem}
For any $N \in \N^*$ the differential equation \eqref{EDO} has a unique 
positive solution $h_N : (0,+\infty) \to (0,+\infty)$. If $N \ge 2$, 
this solution is strictly decreasing and satisfies
\begin{equation}\label{hN1}
  h_N(r) \,=\, 1 + \frac{N{-}1}{2r} + \OO\Bigl(\frac{1}{r^2}\Bigr)
  \qquad \hbox{as} \quad r \to +\infty~,
\end{equation}
and
\begin{equation}\label{hN2}
  h_N(r) \,\sim\, \left\{\begin{array}{ccc}\DS 
  \frac{1}{r\log(1/r)} & \hbox{if} & N = 2~, \\[5mm]
  \DS \frac{N{-}2}{r} & \hbox{if} & N \ge 3~,
  \end{array}\right. \qquad \hbox{as} \quad r \to 0~.
\end{equation}
Moreover, any solution of \eqref{EDO} above $h_N$ blows up in 
finite time, and any solution below $h_N$ cannot stay positive. 
Finally, $h_N$ is given by the explicit formula
\begin{equation}\label{hN3}
  h_N(r) \,=\, \frac{K_{\frac{N}{2}}(r)}{K_{\frac{N}{2}-1}(r)}~, 
  \qquad r > 0~,
\end{equation}
where $K_\nu$ denotes the modified Bessel function as defined in 
\cite[Section~9.6]{AS}. In particular $h_1(r) = 1$ and $h_3(r) = 
1 + 1/r$ for all $r > 0$. 
\end{lem}

We are now able to state the main result of this section\:

\begin{prop}\label{fluxN}
Assume that $N \ge 2$, and let $u(t) = \Phi(t)u_0$ be a trajectory 
for which the energy density $e(x,t)$ satisfies \eqref{e0def}
for some $e_0 > 0$ and $\beta > 0$. Then, for any $R > 0$ and any 
$T > 0$, the integrated energy flux \eqref{intfluxN} satisfies
\begin{equation}\label{FNbd}
   \frac{F(R,T)}{\omega_N R^{N-1}} \,\le\, \sqrt{\beta T e_0}
   ~h_N\Bigl(R \sqrt{\frac{e_0}{\beta T}}\Bigr)~,
\end{equation}
where $h_N$ is given by \eqref{hN3}.
\end{prop}

\noindent{\bf Proof.} Given $T > 0$ and $R_0 > 0$, we shall prove 
inequality \eqref{FNbd} for $R = R_0$. To do that, we consider the 
energy balance equation \eqref{EB} in the spherical shell $\Omega = 
\{x \in \R^N\,|\, R_0 < |x| < R\}$ over the time interval $[0,T]$, 
where $R > R_0$. Using the notation \eqref{intfluxN} we obtain
\begin{equation}\label{FN1}
  F(R,T) \,=\, F(R_0,T) + \int_\Omega \Bigl(e(x,T)-e(x,0)\Bigr)\dd x 
  + \int_0^T \int_\Omega d(x,t)\dd x\dd t~.
\end{equation}
To estimate the right-hand side of \eqref{FN1}, we proceed as 
in the proof of Proposition~\ref{flux1}. We first observe that
\begin{equation}\label{FN2}
  \int_\Omega \Bigl(e(x,T)-e(x,0)\Bigr)\dd x \,\ge\, - e_0 |\Omega| \,=\, 
  - e_0 \frac{\omega_N}{N}(R^N - R_0^N)~.
\end{equation}
Next, applying the Cauchy-Schwarz inequality to \eqref{intfluxN} and 
using \eqref{fbound2}, we obtain the inequality
\[
  |F(R,T)|^2 \,\le\, \omega_N R^{N-1} T \int_0^T \!\!\int_{|x|=R} 
  |f(x,t)|^2 \dd\sigma \dd t 
  \,\le\, \beta T \omega_N R^{N-1} \int_0^T \!\!\int_{|x|=R} 
  d(x,t) \dd\sigma \dd t~,
\]
which implies that
\begin{equation}\label{FN3}
  \int_0^T \!\int_\Omega d(x,t)\dd x\dd t \,=\, \int_0^T \!\int_{R_0}^R
  \int_{|x|=r} d(x,t)\dd \sigma \dd r\dd t \,\ge\, \frac{1}{\beta 
  T \omega_N}\int_{R_0}^R \frac{F(r,T)^2}{r^{N-1}}\dd r~.
\end{equation}
Thus, combining \eqref{FN1} with \eqref{FN2} and \eqref{FN3}, we obtain
\begin{equation}\label{FN4}
  F(R,T) \,\ge\, F(R_0,T) - e_0 \frac{\omega_N}{N}(R^N - R_0^N) +
  \frac{1}{\beta T \omega_N}\int_{R_0}^R \frac{F(r,T)^2}{r^{N-1}}\dd r~, 
  \qquad R > R_0~.
\end{equation}

Now, arguing as in the proof of Proposition~\ref{flux1}, we compare 
$F(R,T)$ to the solution of the differential equation
\begin{equation}\label{FN5}
  \tilde F'(R) \,=\, - e_0 \omega_N R^{N-1} +  \frac{1}{\beta T \omega_N}
  \frac{\tilde F(R)^2}{R^{N-1}}~, \qquad R > R_0~,
\end{equation}
with initial data $\tilde F(R_0) = F(R_0,T)$. To eliminate all parameters
in \eqref{FN5}, we set
\[
  \frac{\tilde F(R)}{\omega_N R^{N-1}} \,=\, \sqrt{\beta T e_0}
  ~H\left(R \sqrt{\frac{e_0}{\beta T}}\right)~, \qquad R > R_0~, 
\]
so that $H$ satisfies the normalized equation \eqref{EDO}. By 
Lemma~\ref{EDOlem}, if 
\[
  \frac{\tilde F(R_0)}{\omega_N R_0^{N-1}} \,\equiv\,
  \frac{F(R_0,T)}{\omega_N R_0^{N-1}} \,>\, \sqrt{\beta T e_0}
  ~h_N\!\left(R_0 \sqrt{\frac{e_0}{\beta T}}\right)~,
\]
then the solution $\tilde F$ of \eqref{FN5} is strictly positive 
for $R > R_0$ and blows up at some finite point $R_* > R_0$. 
But in view of \eqref{FN4}, \eqref{FN5} we also have $\tilde F(R) \le 
F(R,T)$ for all $R \in (R_0,R_*)$, which is impossible since 
$F(R,T)$ is uniformly bounded for all $R \in [R_0,R_*]$. 
Thus we must have 
\[
  \frac{F(R_0,T)}{\omega_N R_0^{N-1}} \,\le\, \sqrt{\beta T e_0}
  ~h_N\!\left(R_0 \sqrt{\frac{e_0}{\beta T}}\right)~,
\]
which is the desired bound. \QED

\begin{rems}\label{fluxrem}\quad\\
{\bf 1.} The bound \eqref{FNbd} also holds for $N=1$ if we set
$\omega_1 = 2$ and $h_1 \equiv 1$. It then asserts that the total
energy entering the segment $[-R,R]$ over the time interval $[0,T]$ is
bounded from above by $2(\beta Te_0)^{1/2}$. This is of course an
immediate consequence of Proposition~\ref{flux1}. 

\noindent{\bf 2.} It is also possible to obtain a lower bound on the
energy flux leaving the ball $B_R$ over the time interval $[0,T]$.
For instance, using the energy balance equation \eqref{EB} with
$\Omega = B_R$ and $[T_1,T_2] = [0,T]$, we easily obtain the estimate
$F(R,T) \ge -\omega_N R^N e_0/N$, which is
uniform in $T$. 

\noindent{\bf 3.} If $N=2$, it follows from \eqref{hN2} that
\begin{equation}\label{Fasym2}
  \sqrt{\beta T e_0}~h_N\Bigl(R \sqrt{\frac{e_0}{\beta T}}\Bigr)
  ~\sim~ \frac{2\beta T}{R\log\Bigl(\frac{\beta T}{e_0 R^2}\Bigr)}
  \qquad \hbox{as}\quad T \to +\infty~.  
\end{equation}
In view of Proposition~\ref{fluxN}, it follows that, for any 
given $R > 0$, the integrated flux $F(R,T)$ can grow at most
sub-linearly (like $T/\log T$) as $T \to +\infty$. This is enough 
to preclude the existence of nontrivial time-periodic solutions in 
two-dimensional extended dissipative systems, using the same 
argument as in Remark~\ref{noper1}, see also \cite{GS}. In contrast, 
if $N \ge 3$, we have
\begin{equation}\label{FasymN}
  \sqrt{\beta T e_0}~h_N\Bigl(R \sqrt{\frac{e_0}{\beta T}}\Bigr)
  ~\sim~ (N-2) \frac{\beta T}{R} \qquad \hbox{as}\quad T \to +\infty~.  
\end{equation}
In that case $F(R,T)$ may grow linearly in time as $T \to \infty$,
which is compatible with the existence of nontrivial time-periodic
orbits (see \cite{GS} for explicit examples). But such solutions must
be spatially localized, because estimate \eqref{FNbd} shows that the
flux per unit area $F(R,T)/(\omega_N R^{N-1})$ decreases like $1/R$ as
the radius $R$ of the sphere increases to infinity.

\noindent{\bf 4.} The right-hand side of \eqref{FNbd} always decreases
to zero as $\beta \to 0$. On the other hand, it is easy to verify that
the asymptotics \eqref{Fasym2} and \eqref{FasymN} also hold for a fixed
$T > 0$ in the limit where $e_0 \to 0$. In particular, if $N \ge 3$,
Proposition~\ref{fluxN} does not preclude the existence of nontrivial
solutions emerging from initial data with zero energy density.

\noindent{\bf 5.} As was mentioned in Remark~\ref{edsrem}.5, when
solving a nonlinear PDE it is often convenient to use a function space
where the energy density is not continuous, but only locally
integrable. Although we do not want to address such technicalities in
the present paper, it is perhaps instructive to see how
Proposition~\ref{fluxN} is modified if we only suppose that the initial
energy is bounded in the {\em uniformly local} sense, namely
\[
  \bar e_0 \,=\, \sup_{x \in \R^N} \int\limits_{|y-x| \le 1} e(y,0)\dd y
  \,<\, \infty~.
\]
For simplicity, we still assume that \eqref{fbound2} holds for some
$\beta > 0$. Arguing as in the proof of Proposition~\ref{fluxN}, 
we fix $R > R_0 \ge 1$ and we consider the energy balance equation
in the spherical shell $\Omega = \{x \in \R^N\,|\, R_0 < |x| < R\}$
over the time interval $[0,T]$. The main difference comes from the
estimate of the energy initially contained in $\Omega$. Using the
definition of $\bar e_0$, we find
\[
  \int_\Omega e(x,0)\dd x \,\le\, c_N \bar e_0\frac{\omega_N}{N}
  (R^N - R_0^N) + d_N\bar e_0\omega_N R_0^{N-1}~,
\]
where $c_N, d_N$ are positive constants related to the optimal 
covering of a (large) ball or sphere in $\R^N$ with balls of 
unit radius \cite{CS}. Inserting this estimate in \eqref{FN2} and
proceeding as before, we can show that $F(R,T) \ge \tilde F(R)$ 
for $R \ge R_0$, where $\tilde F$ is the solution of the ODE
\[
  \tilde F'(R) \,=\, - \bar e_0 c_N \omega_N R^{N-1} +  
  \frac{1}{\beta T \omega_N}\frac{\tilde F(R)^2}{R^{N-1}}~, 
  \qquad R > R_0~,
\]
with initial data $\tilde F(R_0) = F(R_0,T) - \bar e_0 d_N \omega_N
R_0^{N-1}$. This leads finally to the upper bound
\begin{equation}\label{FNbdd}
  \frac{F(R,T)}{\omega_N R^{N-1}} \,\le\, d_N \bar e_0 + 
  \sqrt{c_N \beta T \bar e_0} ~h_N\!\left(R \sqrt{\frac{c_N {\bar e_0}}
  {\beta T}}\right)~, \qquad R \ge 1~,
\end{equation}
which replaces \eqref{FNbd}. Note that the asymptotics as $T \to 
\infty$ are still given by \eqref{Fasym2}, \eqref{FasymN}.  
\end{rems}

\section{Bounds on the Energy Dissipation}\label{secEn}

As is clear from the balance equation \eqref{EB}, a bound on the
amount of energy entering the ball $B_R = \{x \in \R^N \,|\, |x| <
R\}$ over the time interval $[0,T]$ implies an estimate of the energy
dissipated in $B_R$ during the same time, provided the initial energy
in $B_R$ is under control. In this section, we derive various dissipation 
estimates using the flux bounds established in Section~\ref{secFlux}. 
We also show that, for nonequilibrium solutions of extended 
dissipative systems on $\R$ or $\R^2$, energy dissipation must occur 
``almost everywhere'' in space. 

\subsection{Energy dissipation in fixed or increasing domains}

As in Section~\ref{secFlux}, we consider a trajectory $u(t) = \Phi(t)
u_0$ of an extended dissipative system satisfying the uniform bounds
\eqref{e0def} for some $e_0 > 0$ and $\beta > 0$. Given $R > 0$ and 
$T > 0$, we denote by $F(R,T)$ the energy entering the ball $B_R$ 
(through the boundary $\partial B_R$) over the time interval $[0,T]$.
This quantity is defined by \eqref{intfluxN} for $N \ge 2$, and 
if $N = 1$ we set $F(R,T) = F_1(R,T) - F_1(-R,T)$, where $F_1$ is 
given by \eqref{intflux1}. Using the energy balance equation \eqref{EB} 
and Proposition~\ref{fluxN}, we easily obtain
\begin{align}\label{FN6}
  \int_0^T \!\int_{B_R} d(x,t) \dd x \dd t \,&=\, F(R,T) +     
  \int_{B_R} \Bigl(e(x,0)\dd x - e(x,T)\Bigr)\dd x \\ \label{dissip0}
  \,&\le\,  \omega_N R^{N-1} \sqrt{\beta T e_0}~h_N\Bigl(R 
  \sqrt{\frac{e_0}{\beta T}}\Bigr) + \frac{\omega_N}{N}R^{N}e_0~,
\end{align}
where $h_N$ is given by \eqref{hN3}. Equivalently, if $\tilde h_N(r) 
= N h_N(r)/r$, we find
\begin{equation}\label{dissip00}
   \frac{N}{\omega_N R^N} \int_0^T \!\int_{B_R} d(x,t) \dd x \dd t 
   \,\le\, e_0 \left(\tilde h_N \Bigl(R\sqrt{\frac{e_0}{\beta T}}
   \Bigr) + 1\right)~.
\end{equation}
We now investigate a few consequences of the general bound 
\eqref{dissip0} or \eqref{dissip00}. 

First, we fix $R > 0$ and consider the limit where $T \to +\infty$. 
Using the asymptotics \eqref{hN2} for $N \ge 2$ and the fact that
$h_1 = 1$, we obtain the following result. 

\begin{cor}\label{endissip} Under the assumptions of 
Proposition~\ref{flux1} or \ref{fluxN}, the following inequalities 
hold for any $R > 0$. \\[1mm]
1) If $N = 1$, 
\begin{equation}\label{dissip1}
  \limsup_{T \to \infty} \frac{1}{\sqrt{T}} \int_0^T \!\int_{B_R}
  d(x,t)\dd x\dd t \,\le\, 2\sqrt{\beta e_0}~. 
\end{equation}
2) If $N = 2$, 
\begin{equation}\label{dissip2}
  \limsup_{T \to \infty} \frac{\log T}{T} \int_0^T \!\int_{B_R}
  d(x,t)\dd x\dd t \,\le\, 4\pi\beta~. 
\end{equation}
3) If $N\ge 3$, 
\begin{equation}\label{dissip3}
  \limsup_{T \to \infty} \frac{1}{T} \int_0^T \!\int_{B_R}
  d(x,t)\dd x\dd t \,\le\, \beta (N-2)\omega_N R^{N-2}~. 
\end{equation}
\end{cor}

In particular, if $N \le 2$, it follows from \eqref{dissip1},
\eqref{dissip2} that the energy dissipation in any fixed ball
converges to zero ``on average'' as time goes to infinity. Since we
assumed that energy dissipation vanishes only on equilibria of the
system (see assumption~A3 in Definition~\ref{edsdef}), these estimates
will imply that the trajectory $u(t)$ converges ``on average'' 
to the set of equilibria as $t \to \infty$, in a sense that will 
be specified in Section~\ref{secConv}. Observe also that the bounds 
\eqref{dissip1} and \eqref{dissip2} are independent of the 
radius $R$ of the ball, whereas \eqref{dissip2} and \eqref{dissip3}
do not depend on the initial energy density $e_0$. 

It is also instructive to estimate the energy dissipation in a ball
whose radius depends on the observation time $T$. In view of 
\eqref{dissip00}, it is natural to take $R = R_0 \sqrt{T}$ for 
some $R_0 > 0$. We thus find\:

\begin{cor}\label{endissip2} Under the assumptions of 
Proposition~\ref{flux1} or \ref{fluxN}, the following inequality 
holds for any $N \in \N^*$, any $R_0 > 0$, and any $T > 0$\:
\begin{equation}\label{dissip4}
  \frac{N}{\omega_N R_0^N T^{N/2}} \int_0^T \!\int_{B_{R_0\sqrt{T}}}
  d(x,t)\dd x\dd t \,\le\, e_0 \left(\tilde h_N \Bigl(R_0 
  \sqrt{\frac{e_0}{\beta}}\Bigr) + 1\right)~,
\end{equation}
where $\tilde h_N(r) = N h_N(r)/r$. 
\end{cor}

Observe that the volume of the space-time cylinder $B_R \times [0,T]$ 
is $\omega_N N^{-1}R^N T$, hence \eqref{dissip4} implies that the 
energy dissipation rate $d(x,t)$ is very small on average on 
$B_R \times [0,T]$ if $T \gg 1$. This remark will be exploited 
in Section~\ref{secNS}, on a particular example, to prove convergence 
to equilibria uniformly on large domains (whose size increases with
time). 

Finally, in the two-dimensional case, it is also useful to consider
the energy dissipation on a ball whose radius $R$ has a slower 
growth than $T^{1/2}$ as $T \to \infty$. Obvious possibilities are
$R = R_0 T^\gamma$ for $\gamma < 1/2$, or $R = R_0 T^{1/2}/\log(T)$.
This gives the following estimates, which complement \eqref{dissip2}. 

\begin{cor}\label{endissip3} Assume that $N = 2$. Under the assumptions
of Proposition~\ref{fluxN}, the following inequalities hold\:\\[1mm]
1) If $R(T) = R_0 T^\gamma$ for some $R_0 > 0$ and some $\gamma \in 
[0,1/2)$, then
\[
  \limsup_{T \to \infty}\frac{\log(T)}{T} \int_0^T \!\int_{B_{R(T)}}
  d(x,t)\dd x\dd t \,\le\, \frac{4\pi\beta}{1-2\gamma}~.
\]
2)  If $R(T) = R_0 T^{1/2}/\log(T)$ for some $R_0 > 0$, then
\[
  \limsup_{T \to \infty}\frac{\log(\log(T))}{T} \int_0^T \!
  \int_{B_{R(T)}} d(x,t)\dd x\dd t \,\le\, 2\pi\beta~.
\]
\end{cor}

\subsection{Spatial distribution of energy dissipation}

In Section~\ref{secFlux}, we have seen a first way to exploit 
the energy relation \eqref{eb} and the flux bound \eqref{fbound2}
to derive useful information on the dynamics of the system. 
We now consider the problem from a somewhat broader perspective. 
Let $u(t) = \Phi(t)u_0$ be a trajectory of an extended dissipative
system in the sense of Definition~\ref{edsdef}, and suppose that
the energy flux satisfies \eqref{fbound2} for some $\beta > 0$.
This is the case if the function $e \mapsto b(e)$ in assumption~A2
is bounded from above, or in more general situations if the energy 
density $e(x,t)$ is uniformly bounded. Given $R > 0$ and $T > 0$, 
we consider the integrated energy flux $F(R,T)$ defined by 
\eqref{intfluxN}, and we also denote
\[
  E(R,T) \,=\, \int_{B_R} e(x,T)  \dd x~, \qquad 
 \delta E(R,T) \,=\, E(R,T) - E(R,0)~. 
\]
where as usual $B_R = \{x \in \R^N\,|\, |x| < R\}$. With these 
notations, it follows from \eqref{FN1} and \eqref{FN3} that
\begin{equation}\label{Ric1}
  F(R,T) \,\ge\, F(R_0,T) + \delta E(R,T) - \delta E(R_0,T) 
  + \frac{1}{\beta T \omega_N}\int_{R_0}^R \frac{F(r,T)^2}{r^{N-1}}\dd r~,
\end{equation}
for all $R > R_0 > 0$. If we had equality in \eqref{Ric1}, we could
differentiate both sides with respect to $R$ and obtain, after 
omitting the $T$-dependence, the Riccati differential equation
\begin{equation}\label{Ric2}
  \tilde F'(R) \,=\, \delta E'(R) + \frac{\tilde F(R)^2}{\beta T 
  \omega_N R^{N-1}}~, \qquad R > 0~.
\end{equation}
As a matter of fact, if $\tilde F$ is the solution of \eqref{Ric2}
with initial data $\tilde F(R_0) = F(R_0,T) \ge 0$, it follows from
\eqref{Ric1} that $F(R,T) \ge \tilde F(R)$ for $R > R_0$, as long as $
F(R,T) + \tilde F(R) \ge 0$. This comparison principle imposes strong
constraints to the possible solutions of \eqref{Ric1}, because (as we
have already seen) the solutions of the Riccati equation may blow-up
in finite time. Unfortunately, the conditions preventing a blow-up are
not easy to specify in general, because the solutions of \eqref{Ric2}
cannot be written in explicit form. 

One way to proceed is to make simple assumptions on the source term
$\delta E'(R)$ allowing to obtain an explicit solution of
\eqref{Ric2}, to which the solution of \eqref{Ric1} can then be
compared. In Section~\ref{secFlux}, for instance, we assumed that
$\delta E'(R) \ge -\omega_N R^{N-1} e_0$ for some $e_0 > 0$, and we
obtained as a consequence the upper bound \eqref{FNbd}. Here we use
the same strategy to prove that, if $N \le 2$, the energy difference
$\delta E(R,T)$ must be negative for most values of the radius $R >
0$. Given $T > 0$, denote
\begin{equation}\label{Jdef}
  J_T \,=\, \Bigl\{R > 0 \,\Big|\, E(R,T) \ge E(R,0)\Bigr\} 
  \,\subset\, (0,\infty)~.
\end{equation}
The main result of this section is\:

\begin{prop}\label{endiff}
Assume that $u_0 \in X$ is not an equilibrium, and that the 
trajectory $u(t) = \Phi(t)u_0$ satisfies \eqref{fbound2} for 
some $\beta > 0$. Then for any $T > 0$ we have
\begin{equation}\label{Jint}
  \int_1^\infty \frac{\1_{J_T}(r)}{r^{N-1}}\,\dd r \,<\, \infty~,
\end{equation}
where $\1_{J_T}$ is the characteristic function of the set 
$J_T$ defined in \eqref{Jdef}. 
\end{prop}

\begin{rem}\label{endiffrem}
Of course, the conclusion of Proposition~\ref{endiff} is interesting
only if $N \le 2$. If $N = 1$ then \eqref{Jint} simply means that the 
Lebesgue measure of the set $J_T \subset (0,\infty)$ is finite. If
$N = 2$ we have
\[
  \int_1^\infty \frac{\meas(J_T \cap [1,r])}{r^2} \dd r \,=\, 
  \int_1^\infty \frac{\1_{J_T}(r)}{r}\,\dd r \,<\, \infty~,
\]
which implies (roughly speaking) that $\meas(J_T \cap [1,R]) =
o(R/\log(R))$ as $R \to \infty$. In both cases \eqref{Jint} shows that
$J_T$ is a very sparse subset of the half-line $(0,+\infty)$, so that
$E(R,T) < E(0,T)$ for most values of $R > 0$. This considerably
strengthens the results obtained (on a particular example) in
\cite[Section~2]{GS}.
\end{rem}

\noindent{\bf Proof of Proposition~\ref{endiff}.}
Fix $T > 0$. We start from the energy balance equation \eqref{FN6}, 
namely
\begin{equation}\label{Ric3}
  F(R,T) \,=\, \delta E(R,T) + \int_0^T \int_{B_R} d(x,t)\dd x\dd t~,
  \qquad R > 0~. 
\end{equation}
Since $u_0 \in X$ 
is not an equilibrium, assumption~A3 in Definition~\ref{edsdef}
implies that the last term in \eqref{Ric3} is positive when 
$R \ge R_1$, for some (sufficiently large) $R_1 > 0$. If $J_T 
\subset (0,R_1]$, then obviously \eqref{Jint} holds. If this is
not the case, we choose $R_2 \in J_T \cap (R_1,+\infty)$ and 
\eqref{Ric3} then implies that $F(R_2,T) > 0$. Taking the limit 
$R_0 \to 0$ in \eqref{Ric1} we also have
\begin{equation}\label{Ric4}
  F(R,T) \,\ge\, \delta E(R,T) + \frac{1}{\beta T \omega_N}
  \int_0^R \frac{F(r,T)^2}{r^{N-1}}\dd r~, \qquad R > 0~.
\end{equation}
We now define
\[
  \FF(R) \,=\, \frac{1}{(\beta T \omega_N)^2}\int_0^R 
  \frac{F(r,T)^2}{r^{N-1}}\dd r~, \qquad R > 0~.
\]
The function $\FF : (0,\infty) \to \R_+$ is nondecreasing 
and $\FF(R) > 0$ for all $R \ge R_2$. Moreover, using \eqref{Ric4}
and the definition \eqref{Jdef} of $J_T$, we easily find
\[
  \FF'(R) \,\ge\, \1_{J_T}(R)\,\frac{\FF(R)^2}{R^{N-1}}~, \qquad
  R > 0~.
\]
Thus, for all $R > R_2$, we have
\[
  \int_{R_2}^R \frac{\1_{J_T}(r)}{r^{N-1}}\,\dd r \,\le\,  
  \int_{R_2}^R \frac{\FF'(r)}{\FF(r)^2}\dd r \,=\,
  \frac{1}{\FF(R_2)} - \frac{1}{\FF(R)} \,\le\, \frac{1}{\FF(R_2)}~,
\]
and \eqref{Jint} follows. \QED

\begin{rem}\label{strongerA3}
In the proof of Proposition~\ref{endiff}, we used in fact a slightly
stronger version of assumption~A3 in Definition~\ref{edsdef}\: If
$u(t) = \Phi(t)u_0$ is a trajectory for which the energy dissipation
satisfies $d \equiv 0$ on $\R^N \times [0,\epsilon]$ for some
$\epsilon > 0$, then $u_0$ is an equilibrium. With the original
formulation, the conclusion of Proposition~\ref{endiff} would hold
only for sufficiently large $T > 0$. Actually, in all concrete examples 
we are aware of, both versions of assumption~A3 are equivalent.
\end{rem}

\section{Convergence to Equilibria}\label{secConv}

So far we only considered a single trajectory $u(t) = \Phi(t)u_0$ 
of an extended dissipative system, and under appropriate boundedness 
assumptions we established a few estimates on the transfer and
the dissipation of energy. Now, in the spirit of Remark~\ref{noper1}, 
we want to show that these results impose nontrivial restrictions on
the dynamics of the whole system, at least if the space dimension is not
larger than $2$. In particular, we shall use the topology of the 
underlying space $X$ to formulate convergence results, and to study 
the dynamics of the system in a neighborhood of a given point. 

For definiteness, we assume henceforth that our system is {\em
  bounded} and {\em compact} in the sense of Definition~\ref{bcdef}.
As was already mentioned in Remark~\ref{bdrem}, boundedness is 
easy to achieve in the applications by restricting the space $X$ to a
suitable positively invariant subset on which the energy density is
uniformly bounded, and compactness can then be obtained by equipping
$X$ with a sufficiently weak topology. If $u(t) = \Phi(t)u_0$ is a 
trajectory of our system and if $d(x,t)$ is the corresponding energy 
dissipation rate, we denote for all $R > 0$ and all $T > 0$\:
\begin{equation}\label{Ddef}
  D(R,T) \,=\, \int_0^T \Lambda_R(u(t))\dd t~, \qquad 
  \hbox{where} \quad \Lambda_R(u(t)) \,=\, \int_{B_R} d(x,t) \dd x~.
\end{equation}
Here, as usual, $B_R = \{x \in \R^N\,|\, |x| < R\}$. By assumption~A1
in Definition~\ref{edsdef}, we know that $D(R,T)$ depends continuously
on the initial data $u_0$ in the topology of $X$. 

As a consequence of Corollary~\ref{endissip}, we first estimate the
time spent by any trajectory in a neighborhood of a nonequilibrium
point.

\begin{prop}\label{time}
Consider a bounded extended dissipative system on $\R^N$ with 
$N \le 2$. If $\bar u \in X$ is not an equilibrium point, then 
$\bar u$ has a neighborhood $\VV$ in $X$ such that any trajectory
$u(t) = \Phi(t)u_0$ satisfies
\begin{equation}\label{timeN}
  \limsup_{T \to \infty} \frac{\Psi_N(T)}{T} \int_0^T \1_{\VV}(u(t))
  \dd t \,<\, \infty~,
\end{equation}
where $\Psi_1(T) = \sqrt{T}$, $\Psi_2(T) = \log(T)$, and 
$\1_\VV$ denotes the characteristic function of $\VV \subset X$. 
\end{prop}

\noindent{\bf Proof.} We proceed as in \cite[Section~5.1]{GS}. If
$\bar u \in X$ is not an equilibrium point, then assumption~A3 in
Definition~\ref{edsdef} implies that the trajectory $\bar u(t) = 
\Phi(t)\bar u$ satisfies $\bar D(R,T_0) > 0$ for some $R > 0$ and
some $T_0 > 0$, where $\bar D(R,T_0)$ denotes the energy dissipation 
\eqref{Ddef} for the solution $\bar u(t)$. By continuity, there 
exists $\epsilon > 0$ and a neighborhood $\VV$ of $\bar u$ in $X$ 
such that, for any $u_0 \in \VV$, the solution $u(t) = \Phi(t)u_0$ 
satisfies $D(R,T_0) \ge \epsilon > 0$, where $D(R,T_0)$ is given 
by \eqref{Ddef}. 

Now, let $u(t) = \Phi(t)u_0$ be any trajectory of our system. 
Using the notation \eqref{Ddef}, we have for all $T > 0$\:
\begin{align*}
  \frac{1}{T}&\int_0^{T+T_0} \Lambda_R(u(t))\dd t \,\ge\, \frac{1}{T}
  \int_0^T \left(\frac{1}{T_0}\int_t^{t+T_0} \Lambda_R(u(s))\dd s
  \right)\dd t \\
  \,&\ge\, \frac{1}{T} \int_0^T \1_\VV(u(t))\left(\frac{1}{T_0}
  \int_t^{t+T_0} \Lambda_R(u(s))\dd s \right)\dd t
  \,\ge\, \frac{\epsilon}{TT_0} \int_0^T \1_\VV(u(t))\dd t~,
\end{align*}
hence
\[
  \int_0^T \1_\VV(u(t))\dd t \,\le\, \frac{T_0}{\epsilon}
  \int_0^{T+T_0} \!\int_{B_R} d(x,t)\dd x \dd t~, \qquad T > 0~. 
\]
If we multiply both sides by $1/\sqrt{T}$ (if $N = 1$) or $\log T/
\sqrt{T}$ (if $N = 2$) and take the limit $T \to \infty$, we obtain
\eqref{timeN} using Corollary~\ref{endissip}. 
\QED

\begin{rem}\label{weaker}
In \cite{GS}, the following weaker result was obtained for a particular
system\: If $N \le 2$, any nonequilibrium point has a neighborhood 
$\VV$ in $X$ such that any trajectory $u(t)$ satisfies
\[
  \lim_{T \to \infty} \frac{1}{T} \int_0^T \1_{\VV}(u(t)) \dd t \,=\, 0~.
\]
This of course follows from \eqref{timeN}, which gives a much more 
precise estimate of the fraction of time spent by the trajectory
$u(t)$ in the neighborhood $\VV$. 
\end{rem}

Proposition~\ref{time} was obtained without any compactness 
assumption. If we now suppose that the space $X$ is compact, 
we can use \eqref{timeN} to prove that all trajectories converge
in some sense to the set equilibria. Indeed, given any trajectory 
$u(t) = \Phi(t)u_0$, we can define the omega-limit set
\begin{equation}\label{omdef}
  \omega \,=\, \Bigl\{u \in X \,\Big|\, \exists\,t_n \to +\infty 
  \hbox{ such that } u(t_n) \xrightarrow[n \to \infty]{} u 
  \hbox{ in } X\Bigr\} \,\subset\, X~.
\end{equation}
It is known \cite{Ha} that $\omega$ is nonempty, compact, connected,
fully invariant under the semiflow $\Phi(t)$, and that
$\dist_X(u(t),\omega) \to 0$ as $t \to +\infty$. However, our
assumptions {\em do not} imply that $\omega$ is contained in the set
of equilibria. Counter-examples can indeed be constructed even 
for relatively simple systems such as the Allen-Cahn equation in 
one space dimension, see Example~\ref{ERex} below. Motivated 
by the conclusion of Proposition~\ref{time}, we propose the 
following alternative definition\:

\begin{df}\label{omega}
If $u(t) = \Phi(t)u_0$ is a trajectory of a bounded and compact 
extended dissipative system on $\R^N$ with $N \le 2$, we define
\begin{equation}\label{omegadef}
  \bar\omega \,=\, \Bigl\{\bar u \in X \,\Big|\, \limsup_{T \to \infty} 
  \frac{\Psi_N(T)}{T} \int_0^T \1_{\VV}(u(t))\dd t = \infty 
  \hbox{ for all neighborhoods } \VV \hbox{ of }\bar u\Bigr\}~,
\end{equation}
where $\Psi_1(T) = \sqrt{T}$ and $\Psi_2(T) = \log(T)$.  
\end{df}

In other words, $\bar \omega$ is the set of points in all
neighborhoods of which the trajectory $u(t)$ spends a ``substantial
fraction of the total time''. What is exactly meant by ``substantial''
depends on the space dimension $N$, and is specified by the function
$\Psi_N(T)$. It is clear from the definition that $\bar\omega \subset
\omega$, and Proposition~\ref{time} implies that $\bar \omega$ is
contained in the set of equilibria of our system. More properties of
$\bar \omega$ are collected in our final result\:

\begin{prop}\label{barom}
Let $u(t) = \Phi(t)u_0$ be a trajectory of a bounded and compact 
extended dissipative system on $\R^N$ with $N \le 2$. Then the set 
$\bar \omega \subset X$ defined by \eqref{omegadef} is nonempty, 
compact, and contained in the set of equilibria. Moreover, if $\VV$ 
is any neighborhood of $\bar \omega$ in $X$, then
\begin{equation}\label{baromN}
  \limsup_{T \to \infty} \frac{\Psi_N(T)}{T} \int_0^T \1_{\VV^c}(u(t))
  \dd t \,<\, \infty~.
\end{equation}
\end{prop}

\noindent{\bf Proof.} We proceed as in \cite[Section~5.2]{GS}.  We
first observe that, if $\Gamma \subset X$ is compact and does not
intersect $\bar \omega$, then there exists a neighborhood $\VV$ of
$\Gamma$ such that
\begin{equation}\label{Gamma}
  \limsup_{T \to \infty} \frac{\Psi_N(T)}{T} \int_0^T \1_\Gamma(u(t))
  \dd t \,\le\, \limsup_{T \to \infty} \frac{\Psi_N(T)}{T} \int_0^T 
  \1_\VV(u(t))\dd t \,<\, \infty~.
\end{equation}
Indeed, this property holds by definition if $\Gamma = \{u_1\}$ for
some $u_1 \notin \bar\omega$, and the general case follows by a finite
covering argument. Now, if we take for $\Gamma$ the closure of the
trajectory $\{u(t)\,|\, t \ge 0\}$, then $\Gamma$ is compact and
$T^{-1} \int_0^T \1_\Gamma(u(t)) \dd t = 1$ for all $T > 0$, which is
incompatible with \eqref{Gamma}. Thus we must have $\Gamma \cap
\bar\omega \neq \emptyset$, hence in particular $\bar\omega \neq
\emptyset$. Moreover, it is clear from the definition that
$\bar\omega$ is closed in $X$ and contained in $\Gamma$, hence
$\bar\omega$ is compact. On the other hand, Proposition~\ref{time}
precisely means that $\bar\omega$ is contained in the set of
equilibria. Finally, if $\VV$ is any open neighborhood of $\bar\omega$
in $X$, then $\Gamma \cap \VV^c$ is compact and does not intersect
$\bar\omega$, hence by \eqref{Gamma}
\[
  \limsup_{T \to \infty} \frac{\Psi_N(T)}{T} \int_0^T \1_{\VV^c}(u(t))
  \dd t \,=\, \limsup_{T \to \infty} \frac{\Psi_N(T)}{T} \int_0^T 
  \1_{\VV^c \cap \Gamma}(u(t))\dd t \,<\, \infty~,
\]
which proves \eqref{baromN}. \QED

\begin{rem}\label{invariant}
Since $\bar\omega$ consists of equilibria, it is obvious that
$\Phi(t)\bar\omega = \bar\omega$ for all $t \ge 0$. In fact, 
for any relatively compact trajectory of a continuous semiflow
on a metrizable space $X$, one can prove that the set $\bar\omega$
defined by \eqref{omegadef} is nonempty, compact, and fully 
invariant, see \cite[Proposition~5.4]{GS}. These properties 
are therefore independent of the gradient structure. On 
the other hand, the set $\bar\omega$ (unlike $\omega$) 
is not connected in general, as can be seen from Example~\ref{ERex}. 
\end{rem}

\begin{rem}\label{weaker2}
Instead of $\bar\omega$, the following set was defined in \cite{GS}
(for a particular system)\:
\[
  \tilde\omega \,=\, \Bigl\{\bar u \in X \,\Big|\, \limsup_{T \to \infty} 
  \frac{1}{T} \int_0^T \1_{\VV}(u(t))\dd t > 0 \hbox{ for all 
  neighborhoods } \VV \hbox{ of }\bar u\Bigr\}~.
\]
Clearly $\tilde\omega \subset \bar\omega$, hence Proposition~\ref{barom}
implies that $\tilde\omega$ is contained in the set of equilibria, as 
was proved in \cite[Proposition~5.4]{GS}. It is also known that 
$\tilde\omega \neq \emptyset$, which implies that $\bar \omega 
\neq \emptyset$.
\end{rem}

\begin{ex}\label{ERex}
Under the assumptions of Proposition~\ref{barom}, it is not a priori 
obvious that the usual omega-limit set \eqref{omdef} is not 
necessarily contained in the set of equilibria. In this respect, 
the following example is instructive. We consider the one-dimensional 
reaction-diffusion equation \eqref{rd} with the double-well potential 
$V(u) = \frac14(1-u^2)^2$\:
\begin{equation}\label{rGL}
  \partial_t u \,=\, \partial_x^2 u + u - u^3~, \qquad 
  x \in \R~, \quad t \ge 0~.
\end{equation} 
This system has three constant steady states\: $u_0 = 0$ (which is
unstable), and $u_\pm = \pm 1$ (which are stable). In addition, there
is the ``kink'' solution
\begin{equation}\label{psidef}
  \psi(x) \,=\, \tanh(x/\sqrt{2})~, \qquad x \in \R~,
\end{equation}
which connects $u_-$ at $x = -\infty$ to $u_+$ at $x = +\infty$.
It can be shown that $u_\pm$ and the translates of $\pm \psi$ are 
the only stable steady states of \eqref{rGL} in the space of 
bounded solutions. 

Interesting nonequilibrium solutions of \eqref{rGL} can be constructed
by gluing widely separated kinks. For instance, if $a \gg 1$, the
function
\[
  V_a(x) \,=\, \psi(x-a) - \psi(x+a) + 1~, \qquad x \in \R~,
\]
describes the superposition of a kink $\psi$ located near $x = a$ and
an ``anti-kink'' $-\psi$ near $x = -a$. This is not an equilibrium of
\eqref{rGL}, but it can be shown that the solution of \eqref{rGL} with
initial data $V_a$ stays very close to $V_{a(t)}$ for later times,
provided the parameter $a$ evolves according to the exponential law
$\dot a \simeq -c_1\exp(-c_2 a)$, for some $c_1, c_2 > 0$, see e.g.
\cite{CP}.  This approximation property remains valid as long as both
kinks are widely separated, but when they get close to each other they
``annihilate'' and the solution converges uniformly to $1$ as $t \to
+\infty$.

Using these results and a general procedure that can be found 
e.g. in \cite{Ei}, one can show that there exists a unique eternal
solution $u_\psi : \R \times \R \to \R$ of \eqref{rGL} such that
$u_\psi(0,0) = 0$ and 
\begin{equation}\label{upsidef}
  \sup_{x \in \R} | u_\psi(x,t) - V_{a(t)}(x)| \,\xrightarrow[t \to 
  -\infty]{} 0~, \qquad \hbox{where} \qquad a(t) \to +\infty 
  \quad \hbox{as}\quad t \to -\infty~.
\end{equation}
In fact, one has $a(t) \sim c_2^{-1}\log(|t|)$ as $t \to -\infty$.
This solution converges uniformly to $u_+ = 1$ as $t \to +\infty$, and
uniformly on compact sets to $u_- = -1$ as $t \to - \infty$.  If
$\TT_\loc(\R)$ denotes the topology of uniform convergence on compact
sets of $\R$, it follows that $u_\psi(t)$ realizes a {\em heteroclinic
connection} from $u_-$ to $u_+$ through the symmetric collapse of a
pair of kinks coming from infinity. 

Now, using an idea taken from \cite{CE}, we consider the solution 
$u : \R \times \R_+ \to \R$ of \eqref{rGL} with initial data 
$u_0$ satisfying
\begin{equation}\label{u0def}
  u_0(x) \,=\, (-1)^{n+1} \qquad \hbox{if}\quad b_n \le |x| < b_{n+1}~,
\end{equation}
where $(b_n)_{n \in \N}$ is a strictly increasing sequence satisfying
$b_0 = 0$ and $b_{n+1} \gg b_n$ for all $n \in \N$. Under the
evolution of the parabolic equation \eqref{rGL}, the discontinuities
of the initial data are rapidly smeared out, and replaced by smooth
interfaces of the form \eqref{psidef}, the positions of which slowly
move according to exponential law specified above. In particular,
since $b_1 \gg b_0$, the pair of kinks closest to the origin will
hardly feel the presence of the other kinks, and will therefore evolve
in time like the solution of \eqref{rGL} with initial data
$V_{b_0}$. Once the first pair has disappeared, we are essentially
back to the original configuration, with a central pair of kinks that
is now located near $\pm b_1$. This pair evolves on a much slower time
scale, but will eventually come close to the origin and annihilate,
and the same process will continue forever since we started with
infinitely many kinks. Such a coarsening dynamics was studied for
instance in \cite{ER,R}.

These heuristic considerations lead to the following reasonable 
conjecture\:

\begin{conj}\label{coarsening}
Let $u : \R \times \R_+ \to \R$ be the solution of \eqref{rGL}
with initial data \eqref{u0def}. Then the omega-limit set
of the trajectory $(u(t))_{t \ge 0}$ in the topology $\TT_\loc(\R)$ 
is
\[
  \omega \,=\, \{u_+\,,\,u_-\} \cup \{u_\psi(t)\,|\, t \in \R\}
  \cup \{-u_\psi(t)\,|\, t \in \R\}~,
\]
where $u_\pm = \pm 1$ and $u_\psi$ is the eternal solution of 
\eqref{rGL} defined by \eqref{upsidef}. 
\end{conj}

If this conjecture is true, then $\omega$ consists of two equilibria 
$u_\pm$ and two heteroclinc connections between them. Thus $\omega$ 
is a heteroclinic loop, which is not entirely contained in the set of 
equilibria. In contrast, for the same solution, the modified omega-limit 
set introduced in Definition~\ref{omega} satisfies $\bar \omega = 
\{u_+\,,\,u_-\}$, hence is contained in the set of equilibria. 
Note that Proposition~\ref{time} implies that the number of 
annihilations of pairs of kinks that can occur in the time interval
$[0,T]$ is bounded by $C\sqrt{T}$ for large $T$.  
\end{ex}

\section{The Vorticity Equation in an Infinite Cylinder}
\label{secNS}

In this section we analyze in some detail an interesting example which
does not fit exactly into the framework of Definition~\ref{edsdef}, but
can nevertheless be studied using the techniques developed in
Sections~\ref{secFlux} to \ref{secConv}. We consider the incompressible 
Navier-Stokes equation in the infinite cylinder $\Omega = \R \times 
\T$, where $\T = \R/\Z$. Points in $\Omega$ are denoted by  $x = 
(x_1,x_2)$, where $x_1 \in \R$ is the ``horizontal'' and $x_2 \in \T$ 
the ``vertical'' variable. Our system reads
\begin{equation}\label{NS} 
  \partial_t u + (u\cdot\nabla)u \,=\, \Delta u - \nabla p~, \qquad
  \div u \,=\, 0~,
\end{equation}
where $u : \Omega \times \R_+\to \R^2$ denotes the velocity field and 
$p : \Omega \times \R_+ \to \R$ the pressure field. For each $t \ge 0$, 
both quantities $u(x_1,x_2,t)$, $p(x_1,x_2,t)$ are assumed to be bounded 
in $\Omega$ and $1$-periodic with respect to $x_2$. Since $u$ is 
divergence free, we have in particular
\[
  \partial_1 \int_\T u_1(x_1,x_2,t)\dd x_2 \,=\, 
  \int_\T \partial_1 u_1(x_1,x_2,t)\dd x_2 \,=\,
  - \int_\T \partial_2 u_2(x_1,x_2,t)\dd x_2 \,=\, 0~,
\]
hence the vertical average of the horizontal speed, which we denote by
$\langle u_1\rangle$, does not depend on the horizontal variable $x_1$. 
As is explained for instance in \cite{AM}, it then follows from 
\eqref{NS} that $\partial_t \langle u_1\rangle = 0$, so that $\langle 
u_1\rangle$ is a constant which can be set to zero by an appropriate 
Galilean transformation. Thus, without loss of generality, we can assume 
that
\begin{equation}\label{udecomp}
  u(x_1,x_2,t) \,=\, \begin{pmatrix} 0 \\ m(x_1,t)
  \end{pmatrix} + \begin{pmatrix} \hat u_1(x_1,x_2,t) \\ 
  \hat u_2 (x_1,x_2,t)\end{pmatrix}~, \qquad (x_1,x_2) \in 
  \Omega~, \quad t > 0~,
\end{equation}
where $m = \langle u_2 \rangle$. By construction, we then have 
$\langle \hat u_1 \rangle = \langle \hat u_2 \rangle = 0$. 

In addition to \eqref{NS}, we shall study the evolution equation
for the vorticity $\omega = \partial_1 u_2 - \partial_2 u_1$. 
In view of \eqref{udecomp}, we have
\[
  \omega(x_1,x_2,t) \,=\, \partial_1 m(x_1,t) + \hat \omega(x_1,x_2,t)~, 
  \qquad (x_1,x_2) \in \Omega~, \quad t > 0~,
\]
where $\partial_1 m = \langle \omega\rangle$ and $\hat \omega 
= \partial_1 \hat u_2 - \partial_2 \hat u_1$. It is important to 
realize here that, since we want to consider solutions of \eqref{NS}
which do not necessarily decay to zero as $|x_1| \to \infty$, 
the velocity field $u$ is not entirely determined by the vorticity 
$\omega$. More
 precisely, one can show that the oscillating part $\hat u$ of 
the velocity field is given by a Biot-Savart formula\:
\begin{equation}\label{BS}
  \hat u(x_1,x_2,t) \,=\, \int_\R\int_\T \nabla^\perp K(x_1-y_1,x_2-y_2)
  \hat \omega(y_1,y_2,t) \dd y_2 \dd y_1~,
\end{equation}
where $\nabla^\perp = (-\partial_2,\partial_1)$ and 
\begin{equation}\label{Kdef}
  K(x_1,x_2) \,=\, \frac{1}{4\pi}\log\Bigl(2\cosh(2\pi x_1) 
  - 2\cos(2\pi x_2)\Bigr) - \frac{|x_1|}{2}~, \qquad (x_1,x_2) 
  \in \Omega~,
\end{equation}
see \cite{AM}. However, the vertical average $m = \langle u_2\rangle$ 
cannot be completely expressed in terms of the vorticity, and we only
know that $\partial_1 m = \langle \omega\rangle$. The following 
estimates will be useful. 

\begin{lem}\label{BSlem}
There exists a constant $C_1 > 0$ such that, for any $\omega \in 
L^\infty(\Omega)$, the velocity field $\hat u$ defined by \eqref{BS} 
satisfies 
\begin{equation}\label{BS1}
  \|\hat u\|_{L^\infty(\Omega)} \,\le\, C_1 \|\omega\|_{L^\infty(\Omega)}~.
\end{equation}
Moreover, we have $\hat u_1 = -\partial_2 v$ for some $v \in 
L^\infty(\Omega)$, and there exists $C_2 > 0$ such that
\begin{equation}\label{BS2}
  \|v\|_{L^\infty(\Omega)} \,\le\, C_2 \essup_{x_1 \in \R}\left(
  \int_\T \hat\omega(x_1,x_2)^2 \dd x_2\right)^{1/2}  \,\le\, C_2 
  \|\omega\|_{L^\infty(\Omega)}~.
\end{equation}
\end{lem}

\noindent{\bf Proof.} Here and below we denote by $*$ the convolution
on $\Omega$ (considered as an additive group). As is easily verified, 
the function $K$ defined by \eqref{Kdef} satisfies $\nabla K \in 
L^1(\Omega)$, hence
\[
  \|\hat u\|_{L^\infty(\Omega)} \,=\, \|\nabla^\perp K *\hat\omega\|_{
  L^\infty(\Omega)} \,\le\, \|\nabla K\|_{L^1(\Omega)} \|\hat\omega\|_{
  L^\infty(\Omega)} \,\le\, 2\|\nabla K\|_{L^1(\Omega)} \|\omega\|_{
  L^\infty(\Omega)}~.
\]
Similarly, one can check that $K \in L^1(\Omega)$ and
\[
  M_0 \,:=\, \int_\R \,\sup_{x_2 \in \T} |K(x_1,x_2)| \dd x_1 \,<\, \infty~.
\]
Thus, if we define $v = K * \hat \omega$, we have $\hat u_1 = -\partial_2 
v$ by \eqref{BS} and a standard calculation shows that 
\[
  \|v\|_{L^\infty(\Omega)}^2 \,\le\, M_0 \|K\|_{L^1(\Omega)} 
  \essup_{x_1 \in \R} \int_\T \hat\omega(x_1,x_2)^2 \dd x_2~, 
\]
which is the desired result. \QED

\bigskip
Instead of the Navier-Stokes equation \eqref{NS}, we now consider the
evolution system satisfied by the vorticity $\omega$ and the average
speed $m = \langle u_2\rangle$. As in \cite{AM} we obtain
\begin{equation}\label{mom}
  \left\{\begin{array}{rcl}
  \partial_t m + \partial_1 \langle \hat u_1 \hat u_2\rangle 
  \!\! &=& \!\!\partial_1^2 m~, \qquad x_1 \in \R~,
  \\[1mm]
  \partial_t \omega + u\cdot\nabla \omega \!\! &=& \!\!\Delta \omega~,
  \qquad \,(x_1,x_2) \in \Omega~.
  \end{array}\right.
\end{equation}
Here it is understood that $u$ is given by \eqref{udecomp}, where 
$\hat u_1, \hat u_2$ are obtained from $\omega$ via \eqref{BS}. 
Note that system \eqref{mom} is somewhat redundant, because the
horizontal derivative of the first equation is the vertical 
average of the second one, but as is explained above it is not 
possible to get rid completely of the first equation. Given a 
solution of \eqref{mom}, we define for all $x_1 \in \R$ and $t > 0$\:
\begin{align}\nonumber
  e(x_1,t) \,&=\, \frac12 \int_\T \omega(x_1,x_2,t)^2\dd x_2~, \\
  \label{efdNS}
  f(x_1,t) \,&=\, \frac12 \int_\T \Bigl(\partial_1 \omega^2 - u_1 
  \omega^2\Bigr)(x_1,x_2,t)\dd x_2~, \\\nonumber
  d(x_1,t) \,&=\, \int_\T |\nabla\omega(x_1,x_2,t)|^2 \dd x_2~.
\end{align}
In agreement with the general terminology used in this paper, we shall
call $e$ the energy density, $f$ the energy flux, and $d$ the energy
dissipation rate, although the term ``enstrophy'' would certainly be
more appropriate than ``energy'' in the present context. Using
\eqref{mom}, it is easy to verify that the quantities \eqref{efdNS}
satisfy $\partial_t e = \partial_1 f - d$, which is the
one-dimensional version of \eqref{eb}. On the other hand, if $d \equiv
0$, then certainly $\partial_t \omega \equiv 0$ and $\hat \omega =
\omega - \langle \omega\rangle \equiv 0$.  Then $\hat u \equiv 0$ by
\eqref{BS}, and since $\partial_1^2 m =
\partial_1 \langle \omega\rangle \equiv 0$ it follows from \eqref{mom} 
that $\partial_t m \equiv 0$ too. Thus $d \equiv 0$ only for 
equilibria of system \eqref{mom}. Finally, we have the following
estimate for the energy flux\:

\begin{lem}\label{flem}
There exists a constant $C_3 > 0$ such that
\begin{equation}\label{fbdd}
  |f(x_1)|^2 \,\le\, C_3\Bigl(1 + \sup_{y_1 \in \R}e(y_1)\Bigr)  
  e(x_1)d(x_1)~, \qquad \hbox{for all}\quad x_1 \in \R~.   
\end{equation}
\end{lem}

\noindent{\bf Proof.} 
We fix $x_1 \in \R$ and consider both terms in \eqref{efdNS}
separately. First, using the Cauchy-Schwarz inequality, we easily 
find
\[
  \frac12 \Bigl|\int_\T \partial_1 \omega^2\dd x_2\Bigr| \,=\, 
  \Bigl|\int_\T \omega \partial_1 \omega\dd x_2\Bigr| \,\le\,  
  (2 e)^{1/2} d^{1/2}~.
\]
On the other hand, since $u_1 = \hat u_1 = -\partial_2 v$ by 
Lemma~\ref{BSlem}, we have
\[
  \frac12 \int_\T u_1 \omega^2 \dd x_2 \,=\, -\frac12 \int_\T 
  (\partial_2 v) \omega^2 \dd x_2 \,=\, \int_\T v \omega 
  \partial_2 \omega \dd x_2~,
\]
hence using \eqref{BS2} we conclude
\[
  \Bigl|\frac12 \int_\T u_1 \omega^2 \dd x_2\Bigr| \,\le\, 
  \|v\|_{L^\infty(\Omega)}(2 e)^{1/2} d^{1/2} \,\le\, 
  2C_2 \sup_{y_1 \in \R} e(y_1)^{1/2} \,e^{1/2} d^{1/2}~.  
\]
Combining both estimates we obtain \eqref{fbdd}. \QED

\bigskip
The Cauchy problem for Eq.~\eqref{NS} is globally well-posed in the 
Banach space
\[
  X \,=\, \Bigl\{u \in C^0_\bu(\Omega)^2 \,\Big|\, \div u = 0\Bigr\}~,
\]
equipped with the $L^\infty$ norm, see \cite{AM,GMS,ST,Ze}. If $u(t)$
is the solution of \eqref{NS} with initial data $u_0 \in X$, it is
known that $\|u(t)\|_{L^\infty}$ cannot grow faster than $t^{1/2}$ as
$t \to \infty$ (see \cite{AM} and \eqref{mbdd} below), but otherwise we
have no information on the long-time behavior of the solution. In
particular, uniform boundedness is an open problem, which we hope to
address in a future work. Here our goal is to obtain some information
on the associated vorticity $\omega(t)$. Without loss of generality,
we assume that the initial vorticity $\omega_0 = \curl u_0$ is
bounded, and we denote $M = \|\omega_0\|_{L^\infty}$. Since
$\omega(t)$ evolves according to the advection-diffusion equation
\eqref{mom}, the maximum principle implies that $\|\omega(t)
\|_{L^\infty} \le M$ for all $t \ge 0$. It then follows
from \eqref{BS1} that $\|\hat u(t) \|_{L^\infty} \le C_1 M$ for all $t
\ge 0$, so that the oscillating part of the velocity is under
control. On the other hand, if we apply Duhamel's formula to the first
equation in \eqref{mom}, we obtain
\[
  m(t) \,=\, e^{t\partial_x^2}\,m_0 - \int_0^t \partial_1 
  e^{(t-s)\partial_x^2}\,\hat u_1(s)\hat u_2(s)\dd s~, \qquad t > 0~,
\]
where $m_0 = m(0)$ is the vertical average of the vertical initial speed
$u_0$. The uniform bound on $\hat u(t)$ thus implies
\begin{equation}\label{mbdd}
  \|m_1(t)\|_{L^\infty} \,\le\, \|m_0\|_{L^\infty} + \int_0^t
  \frac{\|\hat u(s)\|_{L^\infty}^2}{\sqrt{\pi(t-s)}}\, \dd s
  \,\le\, \|m_0\|_{L^\infty} + \frac{2\sqrt{t}}{\sqrt{\pi}}\,C_1^2 M^2~,
  \qquad t > 0~,
\end{equation}
hence $\|u(t)\|_{L^\infty}$ grows at most like $t^{1/2}$, as already 
announced. 

Under our assumptions, the energy density and the energy flux defined
by \eqref{efdNS} satisfy the following uniform bounds $e(x_1,t) \le
e_0$ and $f(x_1,t)^2 \le \beta d(x_1,t)$, where
\begin{equation}\label{e0beta}
  e_0 \,=\, \frac12\|\omega_0\|_{L^\infty}^2~, \qquad 
  \hbox{and} \qquad \beta \,=\, C_3e_0(1+e_0)~,
\end{equation}
see Lemma~\ref{flem}. Thus we are exactly in position to apply the
results of Sections~\ref{secFlux} and \ref{secEn}. In particular,
using \eqref{dissip0} with $N = 1$, we obtain\:

\begin{prop}\label{omprop1}
If the initial vorticity is bounded, the solution of \eqref{mom}
satisfies, for all $T > 0$ and all $R > 0$, 
\begin{equation}\label{omconv1}
  \int_0^T \int_{-R}^R \int_\T |\nabla \omega(x_1,x_2,t)|^2 \dd x_2
  \dd x_1 \dd t ~\le~ 2 \sqrt{\beta T e_0} + 2R e_0~,
\end{equation}
where $e_0, \beta$ are given by \eqref{e0beta}. 
\end{prop}
  
Proceeding as in Section~\ref{secConv}, one can then use 
\eqref{omconv1} to show that the vorticity $\omega(x_1,x_2,t)$ 
converges uniformly on compact subdomains toward the set of 
equilibria
\[
  \EE \,=\, \Bigl\{w \in C^0_\bu(\Omega)\,\Big|\,
  \nabla w \equiv 0\,,~|w| \le \|\omega_0\|_{L^\infty}\Bigr\}~.
\]
More precisely, adapting Proposition~\ref{barom} to the present
situation, we infer that, if $\VV$ is any neighborhood of $\EE \subset
X$ in the topology of $C^1(\Omega)$, the fraction of the time interval
$[0,T]$ spent by the trajectory $\omega(t)$ outside $\VV$ does not
grow faster than $C T^{1/2}$ as $T \to \infty$. This already gives
valuable information on the solutions of \eqref{mom}, but combining
\eqref{omconv1} with further a priori estimates one can obtain a
stronger and more precise conclusion. In what follows, we assume that
the solution of \eqref{mom} under consideration satisfies
\begin{equation}\label{M1M2}
  \sup_{x_1 \in \R}\int_\T |\partial^2_{12} \omega(x_1,x_2,t)|^2 
  \dd x_2 \,\le\, M_1^2~, \qquad \hbox{and} \qquad
  \sup_{x_1 \in \R} |m(x_1,t)| \,\le\, M_2(1+t)^\beta~, 
\end{equation}
for all $t \ge 0$, where $M_1, M_2 > 0$ and $\beta \in [0,1/2]$.  The
first estimate in \eqref{M1M2} is verifed for all $t \ge 1$ by any
bounded solution of the vorticity equation, due to parabolic
regularization, and the second estimate with $\beta = 1/2$ is just
\eqref{mbdd}. In fact, it is possible to show that \eqref{M1M2} always
holds with $\beta = 1/6$, see \cite{GS2}, but it is still unclear
whether all solutions of \eqref{mom} with bounded initial data satisfy
\eqref{M1M2} with $\beta = 0$.

\begin{prop}\label{omprop2}
Consider a solution of \eqref{mom} with bounded initial data 
such that \eqref{M1M2} holds for some $\beta \in [0,1/2]$. 
If $\beta \le \alpha \le 1/2$, there exists a constant $K_0 > 0$
such that, for any $T \ge 1$, 
\begin{equation}\label{omzero}
  \meas\biggl\{t \in [0,T] \,\bigg|\, \sup_{|x_1| \le T^{(\alpha+2\beta)/3}}
  \,\sup_{x_2 \in \T} |\omega(x_1,x_2,t)| \,\ge\, \frac{K_0}{T^{(\alpha-\beta)/3}}
  \biggr\} \,\le\, K_0 T^{\alpha + 1/2}~.
\end{equation}
\end{prop}

\begin{rem}\label{omzerorem}
Estimate \eqref{omzero} is useful especially when $\beta < \alpha < 1/2$. 
It then shows that the vorticity $\omega(x_1,x_2,t)$ converges 
uniformly to zero on subdomains of size $\OO(t^{(\alpha+2\beta)/3})$,
at a rate comparable to $t^{-(\alpha-\beta)/2}$, except for possible 
excursions whose probability density decays roughly like $t^{\alpha-1/2}$ 
as $t \to \infty$. The fact that $\omega$ converges to zero, 
and not to a nonzero constant $w$, can be understood as follows. 
If $\omega(x_1,x_2,t) = w \neq 0$ on a sufficiently large spatial domain, 
then $\partial_1 m = \langle \omega \rangle = w$ for $x_1$ in a 
large interval, and this is compatible with the assumed upper bound
\eqref{M1M2} only if $w$ is small enough. Thus, in the particular
case of equation \eqref{mom} with bounded initial data for $\omega$ 
and $m$, the omega-limit set ``on average'' as defined in 
Remark~\ref{weaker2} consists of a single point. 
\end{rem}

\noindent{\bf Proof of Proposition~\ref{omprop2}.}
Applying \eqref{omconv1} with $R = \sqrt{T}$, we see that there
exists $C_4 > 0$ such that
\begin{equation}\label{finom1}
  \int_0^T \int_{-\sqrt{T}}^{\sqrt{T}} \int_\T |\nabla \omega(x_1,x_2,t)|^2 
  \dd x_2 \dd x_1 \dd t ~\le~ C_4 \sqrt{T}~,
\end{equation}
for all $T \ge 1$. Given $\alpha \in [0,1/2]$, we define
\[
  J_\alpha(T) \,=\, \left\{t\in [0,T] \,\bigg|\, 
  \int_{-\sqrt{T}}^{\sqrt{T}} \int_\T |\nabla \omega(x_1,x_2,t)|^2
  \dd x_2 \dd x_1 \,\ge\, \frac{1}{T^\alpha}\right\} \,\subset\, 
  [0,T]~.
\]
It follows from \eqref{finom1} that $\meas(J_\alpha(T)) \le C_4 
T^{\alpha + 1/2}$, for all $T \ge 1$. Our goal is to give a 
uniform bound on the vorticity $\omega(x_1,x_2,t)$ on a large 
spatial domain for all $t \in [0,T]\setminus J_\alpha(T)$. 

We observe that $|\omega(x_1,x_2,t)| \le |g(x_1,t)| + 
h(x_1,t)$ for all $(x_1,x_2) \in \Omega$ and all $t \in [0,T]$, 
where 
\[
  g(x_1,t) \,=\, \int_\T \omega(x_1,x_2,t) \dd x_2~, \qquad
  \hbox{and} \qquad h(x_1,t) \,=\, \Bigl(\int_\T
  |\partial_2 \omega(x_1,x_2,t)|^2 \dd x_2\Bigr)^{1/2}~.
\]
We first bound the average $g = \langle \omega \rangle$. If $L 
\le \sqrt{T}$ and $t \in [0,T] \setminus J_\alpha(T)$, we have
\begin{equation}\label{finom2}
  \int_{-L}^L |\partial_1 g(x_1,t)|^2 \dd x_1 \,\le\, \int_{-L}^L 
  \int_\T |\partial_1 \omega(x_1,x_2,t)|^2 \dd x_2 \dd x_1 \,\le\, 
  \frac{1}{T^\alpha}~.
\end{equation}
Furthermore, we know that $g = \partial_1 m$, where $m(x_1,t)$ 
satisfies \eqref{M1M2} for some $\beta \in [0,1/2]$. Using 
\eqref{M1M2}, \eqref{finom2} and Lemma~\ref{lemma1} below, we thus 
obtain
\[
  \sup_{|x_1| \le L} |g(x_1,t)| \,\le\, \frac{C_5 T^\beta}{L}
  + \frac{(2L)^{1/2}}{T^{\alpha/2}}~, \qquad t \in [0,T] \setminus 
  J_\alpha(T)~,
\]
for some $C_5 > 0$. If we now choose $L = T^{(\alpha+2\beta)/3} \le T^{1/2}$, 
we arrive at
\begin{equation}\label{finom3}
  \sup\Bigl\{|g(x_1,t)| \,\Big|\, |x_1| \le T^{(\alpha+2\beta)/3}\Bigr\}
  \,\le\, \frac{C_6}{T^{(\alpha-\beta)/3}}~, \qquad t \in [0,T] 
  \setminus J_\alpha(T)~,
\end{equation}
for some $C_6 > 0$.

On the other hand, we know that $\int_{-L}^L h(x_1,t)^2\dd x_1 \le 
T^{-\alpha}$ when $t \in [0,T] \setminus J_\alpha(T)$. In addition, 
it follows from \eqref{M1M2} that
\[
  |\partial_1 h(x_1,t)| \,\le\, \Bigl(\int_\T |\partial^2_{12} 
  \omega(x_1,x_2,t)|^2 \dd x_2\Bigr)^{1/2} \,\le\, M_1~,
  \qquad x_1 \in \R~, \quad t \in [0,T]~.
\]
Thus Lemma~\ref{lemma2} below implies that
\begin{equation}\label{finom4}
  \sup_{|x_1| \le L} |h(x_1,t)| \,\le\, C \max\Bigl(\frac{M_1^{1/3}}{
  T^{\alpha/3}}\,,\,\frac{1}{L^{1/2}T^{\alpha/2}}\Bigr) \,=\, 
  \frac{C_7}{T^{\alpha/3}}~, \qquad t \in [0,T] \setminus 
  J_\alpha(T)~,
\end{equation}
for some $C_7 > 0$. Combining \eqref{finom3}, \eqref{finom4} we
obtain
\[
  \sup_{|x_1| \le L}\,\sup_{x_2 \in \T} |\omega(x_1,x_2,t)| \,\le\, 
  \sup_{|x_1| \le L} |g(x_1,t)| + \sup_{|x_1| \le L} |h(x_1,t)| 
  \,\le\, \frac{C_6}{T^{(\alpha-\beta)/3}} + \frac{C_7}{T^{\alpha/3}} 
  \,\le\, \frac{C_8}{T^{(\alpha-\beta)/3}}~,
\]
for all $t \in [0,T] \setminus J_\alpha(T)$. Since $L = 
T^{(\alpha+2\beta)/3}$ and $\meas(J_\alpha(T)) \le C_4 T^{\alpha + 1/2}$, 
this gives \eqref{omzero}. \QED

\medskip
Finally, we state and prove two elementary interpolation lemmas which 
were used in the argument above. 

\begin{lem}\label{lemma1} 
Assume that $g \in C^1([0,L])$ satisfies
\[
  \int_0^L g'(x)^2 \dd x \,\le\, \epsilon~, \qquad \hbox{and}\qquad
  \Bigl| \int_0^L g(x)\dd x\Bigr| \,\le\, M~.
\]
Then $\DS \sup_{0 \le x \le L} |g(x)| \,\le\, \frac{M}{L} + (L
\epsilon)^{1/2}$. 
\end{lem}

\noindent{\bf Proof.} We decompose $g(x) = \bar g + h(x)$, where $
\bar g = L^{-1}\int_0^L g(y)\dd y$. Since $h$ has zero mean over 
$[0,L]$, there exists $x_0 \in [0,L]$ such that $h(x_0) = 0$.  
For all $x \in [0,L]$, we thus have
\[
  |h(x)| \,=\, \Bigl|\int_{x_0}^x h'(y)\dd y\Bigr| \,\le\, 
  |x-x_0|^{1/2} \Bigl(\int_{x_0}^x h'(y)^2 \dd y\Bigr)^{1/2}
  \,\le\, (L\epsilon)^{1/2}~.
\]
Since $|\bar g| \le M/L$ by assumption, we obtain the desired 
result. \QED

\begin{lem}\label{lemma2}
Assume that $h \in C^1([0,L])$ satisfies
\[
  \int_0^L h(x)^2 \dd x \,\le\, \epsilon~, \qquad \hbox{and}\qquad
  \sup_{0 \le x \le L} |h'(x)| \,\le\, M~.
\]
Then $\DS\sup_{0 \le x \le L} |h(x)| \,\le\, \max\Bigl((3M\epsilon)^{1/3}\,,
\,(3\epsilon/L)^{1/2}\Bigr)$.
\end{lem}

\noindent{\bf Proof.} If $x_0 \in [0,L]$ is a point where $|h(x)|$
is maximal, we have
\[
  |h(x)| \,\ge\, \|h\|_{L^\infty} - M|x-x_0|~, \qquad x \in [0,L]~,
\]
where $\|h\|_{L^\infty} = \sup\{|h(x)|\,|\, 0 \le x \le L\}$. By 
straightforward calculations, we thus find
\[
  \epsilon \,\ge\, \int_0^L \Bigl(\|h\|_{L^\infty} - M|x-x_0|\Bigr)_+^2
  \dd x \,\ge\, \frac{1}{3}\min\Bigl( \frac{\|h\|_{L^\infty}^3}{M}\,,\,
  \|h\|_{L^\infty}^2 L\Bigr)~.
\]
This gives the desired result. \QED

\section{Appendix\: Proof of Lemma~\ref{EDOlem}}\label{secApp}

In this final section, we study the positive solutions of the ordinary 
differential equation 
\begin{equation}\label{edo}
  h'(r) + \frac{N{-}1}{r} h(r) \,=\, h(r)^2 - 1~, \qquad r > 0~,
\end{equation}
and we prove Lemma~\ref{EDOlem}. All arguments are quite standard,
and are reproduced here for the reader's convenience. Although the
unique positive solution of \eqref{edo} is given by an explicit
formula which can be found using a Cole-Hopf transformation, we 
find it more instructive to prove the first part of Lemma~\ref{EDOlem}, 
including the asymptotics \eqref{hN2} and \eqref{hN3}, without
using this explicit representation, which will be derived only
at the end. We proceed in several steps\: 

\medskip\noindent{\bf 1.} {\em Construction of the stable manifold}. 
The nonautonomous ODE \eqref{edo} has an asymptotic equilibrium 
$h = 1$ at $r = +\infty$, with a one-dimensional {\em stable 
manifold} which contains precisely the solution $h_N$ we are 
looking for. To construct the stable manifold, we set
\begin{equation}\label{gdef}
  h(r) \,=\, 1 + \frac{N{-}1}{2r} + g(r)~, 
\end{equation}
and obtain for $g$ the ODE
\begin{equation}\label{gdiff}
  g'(r) \,=\, 2 g(r) + g(r)^2 - \frac{(N{-}1)(N{-}3)}{4r^2}~.
\end{equation}
As is easily seen, any solution of \eqref{gdiff} that stays bounded
as $r \to +\infty$ satisfies the integral equation
\begin{equation}\label{gint}
  g(r) \,=\, \int_r^\infty e^{2(r-s)}\left(\frac{(N{-}1)(N{-}3)}{4s^2}
  - g(s)^2\right)\dd s~.
\end{equation}
Now, fix $R \in (0,1)$ and take $r_0 > 0$ large enough so that 
$(N{-}1)(N{-}3) \le 4R r_0^2$. It is then straightforward to verify 
that the right-hand side of \eqref{gint} defines a strict contraction
in the closed ball 
\[
  B_{r_0}(R) \,=\, \Bigl\{g \in C^0([r_0,+\infty))\,\Big|\, 
  \sup_{r \ge r_0}|g(r)| \le R \Bigr\}~,
\]
hence has a unique fixed point $g_N \in B_{r_0}(R)$ which, by 
construction, is a solution of \eqref{gdiff} for $r > r_0$. Since
$R$ can be taken arbitrarily small (at the expense of choosing 
$r_0$ large enough), it is clear that $g_N(r) \to 0$ as $r \to 
\infty$. Thus defining
\begin{equation}\label{hNdef}
  h_N(r) \,=\, 1 + \frac{N{-}1}{2r} + g_N(r)~, \qquad r > r_0~,
\end{equation}
we see that $h_N$ satisfies \eqref{edo} and $h_N(r) \to 1$ as $r \to
\infty$. By construction $h_N$ is the unique solution of \eqref{edo}
that stays bounded as $r \to +\infty$.

\begin{rem}\label{trivcase}
If $N = 1$ or $N = 3$, it is clear from \eqref{gint} that $g_N \equiv 0$, 
so that $h_1(r) = 1$ and $h_3(r) = 1 + 1/r$. 
\end{rem}

\noindent{\bf 2.} {\em Asymptotic behavior as $r \to +\infty$}. 
Using \eqref{gint}, we easily find
\[
  g_N(r) \,=\, \frac{(N{-}1)(N{-}3)}{8r^2} + \OO\Bigl(\frac{1}{r^3}
  \Bigr)~, \qquad \hbox{as}\quad r \to +\infty~.
\]
Thus \eqref{hN1} holds, and in view of \eqref{gdiff} we also have
$g_N'(r) = \OO(r^{-3})$ as $r \to +\infty$, so that
\[
  h_N'(r) \,=\, -\frac{N{-}1}{2r^2} + \OO\Bigl(\frac{1}{r^3}\Bigr)~, 
  \qquad \hbox{as}\quad r \to +\infty~.
\]
If $N \ge 2$, this shows that $h_N'(r) < 0$ for $r > 0$ sufficiently large. 

\medskip\noindent{\bf 3.} {\em Global monotonicity}. We assume from 
now on that $N \ge 2$. Solving \eqref{edo} backwards, we construct
(for some $r_* \ge 0$) a maximal solution $h_N : (r_*,+\infty) \to \R$ 
which coincides with \eqref{hNdef} for $r > r_0$. We claim that 
$h_N'(r) < 0$ for all $r > r_*$. Indeed, assume on the contrary that
there exists $\bar r > r_*$ such that $h_n'(\bar r) = 0$ and 
$h_N'(r) < 0$ for all $r > \bar r$. Then
\[
  h_N''(\bar r) \,=\, \Bigl(2h_N(\bar r) - \frac{N{-}1}{\bar r}
  \Bigr)h_N'(\bar r) + \frac{N{-}1}{\bar r^2}h_N(\bar r) \,>\, 0~,
\]
because $h_N'(\bar r) = 0$ and $h_N(\bar r) > 1$. This implies that
$h_N'(r) > 0$ for $r > \bar r$ close enough to $\bar r$, in 
contradiction with the definition of $\bar r$. Thus $h_N'(r) > 0$ 
for all $r > r_*$, and using \eqref{edo} we deduce
\begin{equation}\label{hNbdd}
  1 \,<\, h_N(r) \,<\, \frac{N{-}1}{2r} + \sqrt{1 + \frac{(N{-}1)^2}{
  4r^2}} \qquad \hbox{for all}\quad r > r_*~.
\end{equation}
This estimate shows in particular that $h_N$ cannot blow up at 
a finite point $r > 0$, hence we necessarily have $r_* = 0$
and $h_N'(r) < 0$ for all $r > 0$. 

\medskip\noindent{\bf 4.} {\em Asymptotic behavior as $r \to 0$}.
Setting $f_N = 1/h_N$ we obtain the ODE
\[
  f_N'(r) \,=\, \frac{N{-}1}{r} f_N(r) + f_N(r)^2 - 1~, 
  \qquad r > 0~,
\]
which is very similar to \eqref{edo}. In particular, we have
\begin{equation}\label{fNint}
  f_N(r) \,=\, r^{N-1}f_N(1) + \int_r^1 \Bigl(\frac{r}{s}\Bigr)^{N-1}
  (1 - f_N(s)^2)\dd s~, \qquad 0 < r < 1~.
\end{equation}
A direct study of \eqref{fNint} gives the following asymptotic 
expansion as $r \to 0$\:
\[
  f_2(r) \,=\, r\log\frac{1}{r} + Cr + \OO\Bigl(r^3\Bigl(\log
  \frac{1}{r}\Bigr)^2\Bigr) \quad \hbox{for some } C \in \R~,
\]
whereas $f_3(r) = r/(1+r)$ and
\[
  f_4(r) \,=\, \frac{r}{2} + \OO\Bigl(r^3\log\frac{1}{r}\Bigr)~,
  \qquad f_N(r) \,=\, \frac{r}{N{-}2} + \OO(r^3) \quad \hbox{if } 
  N \ge 5~. 
\]
Since $h_N = 1/f_N$, this proves \eqref{hN2}. 

\medskip\noindent{\bf 5.} {\em Uniqueness and threshold behavior}.
We first study the solutions of \eqref{edo} that lie above $h_N$. 
Assume that $h$ is a solution of \eqref{edo} such that 
$h(r_1) > h_N(r_1)$ for some $r_1 > 0$. In particular, we have 
that $h(r) > h_N(r) \ge 1$ for all $r > r_1$. If $h'(r) \le 0$ for 
all $r \ge r_1$, then $h(r)$ converges to some limit $\ell \ge 1$ 
as $r \to +\infty$, and using \eqref{edo} we easily see that 
$\ell^2 = 1$, hence $\ell = 1$. This implies that the function 
$g(r)$ defined by \eqref{gdef} is small for sufficiently large
$r > 0$ and satisfies the integral equation \eqref{gint}, hence
coincides with $g_N(r)$, which is of course impossible since 
$g(r_1) > g_N(r_1)$. Thus there must exist $r_2 \ge r_1$ such 
that $h(r_2) > h_N(r_2)$ and $h'(r_2) > 0$. If we now choose $r_3 
> r_2$ so that $h'(r) > 0$ for all $r \in [r_2,r_3]$, we have on 
that interval $h(r) > (N{-}1)/(2r)$, hence 
\[
  h''(r) \,=\, h'(r)\left(2h(r) - \frac{N{-}1}{r}\right)
  + \frac{N{-}1}{r^2}h(r) \,>\, 0~, \qquad r_2 \le r \le r_3~.
\]
This argument shows that $h(r)$ is convex for $r \ge r_2$, and 
blows up at some finite point $r^* > r_2$. Indeed, if $h(r)$ 
was defined for all $r > r_2$, the convexity would imply that 
$h(r) \to +\infty$ as $r \to +\infty$, so that $h$ would satisfy 
the differential inequality $h'(r) \ge \frac12 h(r)^2$ for all 
sufficiently large $r$, which is impossible because this inequality 
has no global positive solutions.

We next consider solutions of \eqref{edo} that lie below $h_N$. 
Assume that $h$ is a solution of \eqref{edo} such that 
$0 < h(r_1) < h_N(r_1)$ for some $r_1 > 0$, so that $h(r) < h_N(r)$
for all $r \ge r_1$. If $h(r) \ge 0$ for all $r \ge r_1$, we 
have $h'(r) < 0$ for all $r \ge r_1$, hence $h(r)$ converges to 
some limit $\ell \in [0,1]$ as $r \to +\infty$. But the same 
arguments as above show that $\ell^2 = 1$ and $\ell \neq 1$, 
which is a contradiction. So the solution $h(r)$ must necessarily 
change sign for $r > r_1$. 

It follows in particular that $h_N$ is the unique positive 
solution of \eqref{edo} that is defined for all $r > 0$. 

\medskip\noindent{\bf 6.} {\em Explicit representation}.
Let $u(r) = \exp(-\int_1^r h_N(s)\dd s)$. Then $u$ solves the 
second order {\em linear} ODE
\[
  u''(r) + \frac{N-1}{r}u'(r) \,=\, u(r)~, \qquad r > 0~,
\]
and $u(r)$ decays exponentially to zero as $r \to +\infty$. 
Setting $\nu = 1 - N/2$ and $u(r) = r^\nu v(r)$, we obtain for
$v$ the differential equation
\[
  r^2 v''(r) + r v'(r) - (r^2 + \nu^2) v(r) \,=\, 0~, \qquad r > 0~,
\]
which defines the modified Bessel functions, see \cite[Eq.~9.6.1]{AS}. 
Since $v(r)$ decays exponentially as $r \to +\infty$, we must 
have $v(r) = C K_\nu(r)$ for some $C > 0$, see \cite[Section~9.7]{AS}. 
Thus $u(r) = C r^\nu K_\nu(r)$, and using \cite[Eq.~9.6.28]{AS}
we also find $u'(r) = -Cr^\nu K_{\nu-1}(r)$. Since $K_{-\nu}(r) = 
K_\nu(r)$ by \cite[Eq.~9.6.6]{AS}, we conclude that
\[
  h_N(r) \,=\, -\frac{u'(r)}{u(r)} \,=\, \frac{K_{\nu-1}(r)}{
  K_\nu(r)} \,=\, \frac{K_{\frac{N}{2}}(r)}{K_{\frac{N}{2}-1}(r)}~, 
  \qquad r > 0~,
\]
which proves \eqref{hN3}. The proof of Lemma~\ref{EDOlem} is now 
complete. \QED

\bigskip\noindent{\bf Acknowledgements.} Part of this work was done
when the second author visited Institut Fourier at University
of Grenoble, whose hospitality is gratefully acknowledged. The
authors thank P. Pol\'a\v cik and A. Scheel for fruitful discussions. 
Th.G. was supported in part by the ANR projet PREFERED of the 
French Ministry of Research, and S.S. by the grant No 037-0372791-2803 
of the Croatian Ministry of Science.

\end{document}